\documentclass[12pt]{amsart}

\usepackage[arrow,dvips,matrix,arrow,ps,color,line,curve,frame]{xy}

\usepackage{amssymb}

\advance\textheight1cm

\let\oldlabel=\label
\def\prellabel{\marginparsep=1em
    \def\label##1{\oldlabel{##1}\ifmmode\else\ifinner\else
         \marginpar{{\footnotesize\ \\ \tt
                    ##1}}\fi\fi}}
\advance\headheight1.15pt
\let\:\colon
\let\epsilon\varepsilon

\let\phi=\varphi
\let\theta=\vartheta

\let\Bbb=\mathbb

\def\opn#1#2{\def#1{\operatorname{#2}}}
\opn\gp{gp} \opn\Max{Max} \opn\Ker{Ker} \opn\Coker{Coker}
\opn\Ext{Ext} \opn\conv{conv} \opn\chara{char} \opn\n{n} \opn\h{h}
\opn\GL{GuL} \opn\SL{SL} \opn\sn{sn} \opn\inte{int} \opn\End{End}
\opn\rank{rank} \opn\aff{aff} \opn\Spec{Spec} \opn\Proj{Proj}
\opn\gr{gr} \opn\QF{QF} \opn\I{Im} \opn\Hom{Hom} \opn\Aut{Aut}
\opn\W{Witt} \opn\w{w} \opn\inte{int} \opn\pyr{pyr} \opn\l{l}
\opn\r{r} \opn\const{const} \opn\M{\bf M} \opn\Mf{\bf \M_f}
\opn\V{Ver} \opn\SS{\mathfrak S} \opn\TT{\mathfrak T}
\opn\UU{\mathfrak U} \opn\Trf{trsfm} \opn\ini{in}

\def\cc{{\mathfrak c}}

\def\ZZ{{\Bbb Z}}
\def\RR{{\Bbb R}}

\def\NN{{\Bbb N}}
\def\QQ{{\Bbb Q}}

\def\Bb{{\mathcal B}}%
\def\Q{{\Box\kern1pt}}%

\def\k{{\bf k}}

\def\1{^{-1}}

\newtheorem{lemma}{Lemma}[section]

\newtheorem{theorem}[lemma]{Theorem}
\newtheorem{proposition}[lemma]{Proposition}

\theoremstyle{definition}
\newtheorem{definition}[lemma]{Definition}
\newtheorem{remark}[lemma]{Remark}

\newtheorem{conjecture}[lemma]{Conjecture}

\begin{document}

\title[The nilpotence conjecture]
{The nilpotence conjecture in\\
$K$-theory of toric varieties}

\author{Joseph Gubeladze}

\thanks{Supported by Mathematical Sciences Research Institute,
INTAS grant 99-00817, and TMR grant ERB FMRX CT-97-0107}


\address{Department of Mathematics, San Francisco
State University, San Francisco, CA 94132, USA}

\email{soso@math.sfsu.edu}

\maketitle

\hfill \emph{\small To Richard G. Swan}
%\hfill{\emph{\small the occasion of his 70th birthday}

\begin{abstract}
It is shown that all nontrivial elements in higher $K$-groups of
toric varieties are annihilated by iterations of the natural
Frobenius type endomorphisms. This is a higher analog of the
triviality of vector bundles on affine toric varieties.
\end{abstract}

\section{Introduction}\label{Mnres}

\subsection{\emph{The conjecture}}
Quillen has shown in \cite{Q1} that a regular ring $R$ and its
polynomial extensions $R[T_1,\ldots,T_n]$ ($n\in\NN$) have the
same higher $K$-groups:
\begin{equation}\label{QUILLEN}
K_p(R)=K_p(R[T_1,\ldots,T_n]),\quad p\ge0.
\end{equation}
This extends classical results of Grothendieck \cite{Se} ($p=0$)
and Bass-Heller-Swan \cite{BHS} ($p=1$). The ring
$R[T_1,\ldots,T_n]$ can be thought of as the monoid $R$-algebra of
the free commutative monoid $\ZZ_+^n$.

Now consider an arbitrary convex cone $C\subset\RR^n$ (i.~e.
$C\subset\RR^n$ is a subset such that $\lambda x+\mu y\in C$
whenever $x,y\in C$ and $\lambda,\mu\in\RR_+$), containing no pair
$\{x,-x\}$, $x\not=0$. In this paper we address the question: what
can be said about the relationship between the higher $K$-groups
of $R$ and the monoid ring $R[C\cap\ZZ^n]$? Observe that the
latter is a regular ring if and only if $R[C\cap\ZZ^n]$ is a
polynomial ring. Our main result (Theorem \ref{yes}) answers this
question when $R$ is a field of characteristic 0.

The  \emph{nilpotence conjecture} in $K$-theory of toric
varieties, treated in our previous papers, asserts the following:

\begin{conjecture}\label{conj}
Let $R$ be a (commutative) regular ring, $M$ be an arbitrary
commutative, cancellative, torsion free monoid without nontrivial
units, and $p$ be a nonnegative integer. Then for every sequence
${\bf c}=(c_1,c_2,\dots)$ of natural numbers $\geq2$ and every
element $x\in K_p(R[M])$ there exists an index $j_0\in\NN$ such
that $(c_1\cdots c_j)_*(x)\in K_p(R)$ for all $j>j_0$.
\end{conjecture}

Here $R[M]$ is the monoid $R$-algebra of $M$ and for a natural
number $c$ the endomorphism of $K_p(R[M])$, induced by the
$R$-algebra endomorphism $R[M]\to R[M]$, $m\mapsto m^c$, $m\in M$
is denoted by $c_*$, writing the monoid operation
multiplicatively. (Recall, a commutative monoid $M$ is called
\emph{cancellative} and \emph{torsion free} if it embeds into a
torsion free abelian group.)

The property of the action of $\NN$ on $K_p(R[M])$, mentioned in
Conjecture \ref{conj}, will be called \emph{the nilpotence
property}.

The following is a reformulation of Conjecture \ref{conj} in the
typical case when $\bf c$ is a constant sequence:

\

\noindent{\bf Conjecture.} Let $R$ and $p$ be as above and $c$ be
a natural number $\ge2$. Assume $C\subset\RR^n$ ($n\in\NN$) is a
convex cone, containing no affine line. Then
$$
K_p(R)=K_p\left(R\left[C\cap\left(\ZZ\left[1/c\right]
\right)^n\right]\right).
$$
(Here $\ZZ[1/c]$ is the additive group of the localization of
$\ZZ$ at $c$.)

\

The motivation behind Conjecture \ref{conj} can be described by
the diagram of relationships:
$$ 
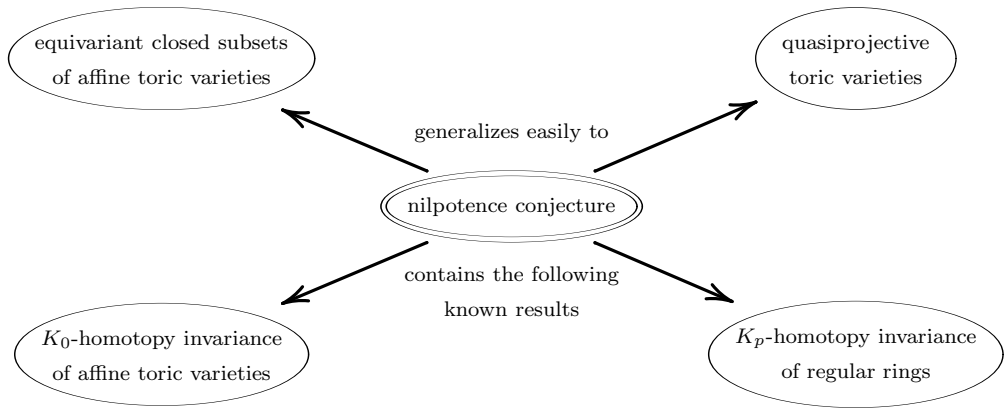
\xymatrix@!R=4pt@M10pt{ \txt{\tiny equivariant closed
subsets\\\tiny of affine toric varieties}\drop\frm{e}%<30pt>{-}
&&\txt{\tiny quasiprojective\\\tiny toric
varieties}\drop\frm{e}%<30pt>{-}
\\
&\txt{\tiny generalizes easily to}\\
&\txt{\tiny nilpotence conjecture}\drop\frm{ee}
\ar@*{[|(3)]}[luu]\ar@*{[|(3)]}[ldd]\ar@*{[|(3)]}[ruu]
\ar@*{[|(3)]}[rdd]\\
&\txt{\\\tiny contains the following\\\tiny known results}\\
\txt{\tiny $K_0$-homotopy invariance\\\tiny of affine toric
varieties}\drop\frm{e}%<30pt>{-}
&&\txt{\tiny $K_p$-homotopy invariance\\\tiny of regular
rings}\drop\frm{e} }
$$
which we are going to describe now.

\noindent$\bullet$ When $p=0$ and $R$ is a PID the nilpotence
conjecture follows from \cite{Gu1}, where we proved the stronger
result on the triviality of algebraic vector bundles on an
arbitrary (normal) affine toric variety. There our starting point
had been Quillen's solution \cite{Q2} to Serre's problem on
projective modules over polynomial rings (solved also by Suslin
\cite{Su1}). That the case of higher dimensional coefficient rings
also follows from \cite{Gu1} was observed by Swan \cite[Corollary
1.4]{Sw2}. Actually, the south-west arrow in our diagram says that
for $p=0$ the nilpotence conjecture is equivalent to the formally
stronger equalities $K_0(R)=K_0(R[M])$ where $M$ runs through all
\emph{normal} monoids\footnote{The monoid terminology will be
introduced in Section \ref{preliminaries}.}. This is proved in
\cite[Proposition 3.5(a,b)]{Gu6}.

\noindent$\bullet$ The case $p=1$ was done in \cite{Gu2} and,
shortly thereafter, Mushkudiani proved the case $p=2$ \cite{M}.
More precisely, what was shown in \cite{Gu2} and \cite{M} is the
aforementioned typical case of the nilpotence conjecture for $K_1$
and $K_2$. The point of departure for us in \cite{Gu2} had been
Suslin's work on the structure of the special linear group of a
polynomial ring \cite{Su2} -- a $K_1$-analog of Serre's problem.

\noindent$\bullet$ While the results for $p\le2$ serve as a strong
evidence for Conjecture \ref{conj}, its non-triviality is
emphasized by the fact that the naive higher analog of \cite{Gu1}
(i.~e. equalities of type $K_p(R)=K_p(R[M])$) is not possible for
essentially any finitely generated monoids $M$ except the
classical case $M=\ZZ^n_+$. First examples of nontrivial elements
in $K_1(R[M])$ for certain rank 2 normal monoids $M$ are due to
Srinivas \cite{Sr}. Later, in \cite{Gu4}, we proved that the
groups $K_1(R[M])/K_1(R)$ are never trivial for the monoid rings
of finitely generated \emph{simplicial} monoids $M$ except when
$M=\ZZ^n_+$. The reason for such a dramatic failure of the direct
higher analog of \cite{Gu1}, as explained in \cite{Gu3}, is the
lack of the excision property for the Bass functor $K_1$ -- a
phenomenon observed by Swan \cite{Sw1}. The examples of Swan have
been studied in detail by Dennis and Krusemeyer \cite{DK}. The
latter work plays a crucial r\^ole in \cite{Gu3} where we
explicitly construct nontrivial elements in $K_1$.

\noindent$\bullet$ The case of simplicial monoids of any rank
(this includes all finitely generated rank 2 monoids without
nontrivial units) was done in \cite{Gu3}. That result is based on
Suslin-Wodzicki's solution \cite{SuW} to the excision problem in
algebraic $K$-theory and, in particular, explains how the
multiplicative action of $\NN$ on $K_p(R[M])$ compensates the lack
of the excision property. On the other hand there is a big
difference between the simplicial and general cases. A further
support for the nilpotence conjecture, obtained in \cite{Gu3}, is
the implication indicated by the north-west arrow in the diagram
above (see also \cite[Proposition 3.5(c)]{Gu6}). This is a higher
analog of Vorst's result \cite{V2} on Serre's problem for discrete
Hodge algebras.

\noindent$\bullet$ Further evidence that Conjecture \ref{conj} is
a natural higher analog of the triviality of vector bundles on
affine toric varieties is the north-east arrow in our diagram: it
is proved in \cite[Proposition 4.7]{Gu5} that the nilpotence of
$K$-theory of affine toric varieties extends to all
quasiprojective toric varieties.

\noindent$\bullet$ The mentioned generalization to quasiprojective
toric varieties makes it clear that the nilpotence property is
essentially \emph{all} that survives from the homotopy invariance
in higher $K$-theory of singular toric varieties. In fact, we have
shown in \cite{Gu7} that even for \emph{projective} simplicial
toric varieties the Grothendieck group of vector bundles can be
`immeasurably' larger than the Grothendieck group of coherent
sheaves, contrary to the previously existing conjecture on the
rational isomorphism between the two groups \cite{BV},
\cite[\S7.12]{Cox}.

\noindent$\bullet$ Finally, Conjecture \ref{conj} actually
contains Quillen's fundamental equality (\ref{QUILLEN}) on which
all of the theory is based -- the south-east arrow in our diagram.
This follows immediately from transfer maps in $K$-theory
\cite[\S1]{Gu6}.

\subsection{\emph{Main result and its proof}}
Until now no example of a coefficient ring was known for which the
nilpotence conjecture is true for all $K$-theoretical indices $p$
and all monoids $M$. In this paper we prove the nilpotence
conjecture for all fields of characteristic 0:

\begin{theorem}\label{yes}
Let $\k$ be a field of characteristic 0, $M$ be an arbitrary
commutative, cancellative, torsion free monoid without nontrivial
units, $p$ be a natural number, and ${\bf c}=(c_1,c_2,\dots)$ be a
sequence of natural numbers $\geq2$. Then for every element $x\in
K_p(\k[M])$ there exists an index $j_0\in\NN$ such that
$(c_1\cdots c_j)_*(x)\in K_p(\k)$ whenever $j>j_0$.
\end{theorem}

The decisive step towards Theorem \ref{yes} has been our previous
work \cite{Gu6} which, using Thomason's Mayer-Vietoris sequence
for singular varieties \cite{TT} and virtually all of our previous
results `almost proves' the nilpotence conjecture when the
coefficient ring is a field. This `almost proves' is made precise
in Section \ref{almost} below. One of the main features of
\cite{Gu6} is the passage to auxiliary \emph{non-affine} toric
varieties (with two affine toric charts) and the endomorphism
rings of certain rank 2 vector bundles on them. This explains the
heavy use of $2\times2$ matrix rings in this paper.

That for a restricted class of nonsimplicial toric
varieties\footnote{including, for instance,
$k[T_1,T_2,T_3,T_4]/(T_1T_2-T_3T_4)$ -- the coordinate ring of the
cone over the Segre embedding ${\Bbb P}^1\times{\Bbb P}^1\to{\Bbb
P}^3$.} over number fields\footnote{In view of Section
\ref{essentially} below the generalization to arbitrary fields of
characteristic 0 is immediate.} the aforementioned `almost proof'
can be converted into an actual proof was also shown in
\cite{Gu6}, based on the Bloch-Stienstra-Weibel action of the big
Witt vectors on the nil-$K$-theory (\cite{Bl}\cite{St}\cite{W}) --
an idea suggested by Burt Totaro, and Goodwillie's rational
isomorphism between relative $K$-groups and relative cyclic
homology groups for a nilpotent ideal \cite{Go} -- a technique of
crucial importance also for this paper. The mentioned action of
the big Witt vectors led us in \cite{Gu6} to a question on certain
filtration on certain relative $K$-groups, which, if answered in
the positive, completes the proof of Conjecture \ref{conj} when
the coefficient ring is a characteristic 0 field \cite[Question
10.2]{Gu6}.

The main part of the present paper is focused on the arithmetic
case and only at the very end (Section \ref{essentially}) we
derive the result for fields of characteristic 0. The reason why
we consider first the arithmetic case is that both in Goodwillie's
and Cortinas' theorems cyclic homology groups are those of
$\ZZ$-algebras, see Remark \ref{WhyQ} for a more detailed
explanation. What is achieved in this work is a finiteness result
on the aforementioned filtration on relative $K$-groups when the
coefficient ring is a number field, and this suffices to complete
the proof of Conjecture \ref{conj} in the arithmetic case.  Such a
progress was made possible by Corti\~nas' recent proof \cite{Cor}
of the Geller-Weibel \emph{KABI conjecture} \cite{GW}, making it
possible to convert our relative $K$-groups into cyclic and,
eventually, Hochschild homology groups whose ranks are amenable to
control.

The general case of coefficient fields is derived by a combination
of Vorst's localization technique \cite{V1}, Galois descent for
rational $K$-theory of rings due to Thomason \cite{T}\cite{TT},
and the aforementioned Bloch-Stienstra-Weibel action of big Witt
vectors.

\subsection{\emph{Future directions}}\label{future}
Here are several possible ways of extending Theorem \ref{yes}.

\noindent$\bullet$ In Section \ref{essentially} we show that the
class of coefficient rings for which Conjecture \ref{conj} is true
is closed under taking polynomial extensions and localizations.
This provides a strong evidence that the result extends to
\emph{all} regular coefficient rings containing $\QQ$. In more
detail, let $n$ be a natural number and $\k$ be a perfect field
(of an arbitrary characteristic). Lindel, in his proof of the
\emph{Bass-Quillen conjecture} in the geometric case \cite{Li},
showed the following implication: the fact that projective
$A[\ZZ_+^n]$-modules are extended from $A$ for the rings of type
$A=S^{-1}\k[\ZZ_+^d]$ ($d\in\NN$ and $S\subset\k[\ZZ_+^d]$)
implies that projective $B[\ZZ_+^n]$-modules are extended from $B$
for all regular rings $B$ essentially of finite type over $\k$.
The proof in \cite{Li} is by induction on $\dim B$, the case $\dim
B=0$ being corresponding to the Quillen-Suslin Theorem. It is very
natural to expect that Lindel's technique of \emph{\'etale
neighborhoods} admits a higher $K$-theoretic analogue which would
extend Theorem \ref{yes} to the coefficient rings which are
regular and essentially of finite type over a field of
characteristic 0. But then, by Popescu's approximation theorem
\cite{P,Sw3}, such a result would automatically extend to all
regular rings containing $\QQ$.

\noindent$\bullet$ Let $\k$ be a field, $M$ be a finitely
generated monoid without nontrivial invertible elements, and $p$
be a nonnegative integer. One could ask whether there is a uniform
bound $\gamma$ such that for any natural number $c>\gamma$ the
endomorphism $c_*$ of $K_p(\k[M])/K_p(\k)$ is the zero map. It is
not difficult to show that such $\gamma$ exists for the functor
$K_0$. If $M$ is normal then one can take $\gamma=0$ (meaning
$K_0(\k[M])/K_0(\k)=0$). Moreover, the explicit nontrivial
elements in $K_1(\k[M])/K_1(\k)$ derived in \cite{Gu4} trivialize
by all endomorphisms $c_*$, $c>1$. We do not know whether all of
the group $K_1(\k[M])/K_1(\k[M])$ trivializes by $c_*$, $c>1$,
though. The question on the existence of such a uniform bound
$\gamma$ seems to be very interesting. But currently it is well
beyond the capacities of the techniques developed so far. A
possible approach may be related to the aforementioned action of
the big Witt vectors $\W(R)$ which provides a framework for
treating the operator $c_*$ as an analogue of the Verschiebung
operator, see Remark \ref{refereetotaro} for details.

\noindent$\bullet$ Other possible higher $K$-theoretic analogs of
\cite{Gu1}, mentioned in \cite[Conjecture 2.4]{Gu5}, concern
stabilizations in higher $K$-groups of monoid rings. The starting
observation here, a sort of homotopy invariance, is the nontrivial
fact that a ring and its polynomial extensions have essentially
the same stabilization behavior.

\subsection{Organization of the paper.}\label{organization} In Section
\ref{preliminaries} we recall the polyhedral geometry, relevant to
\emph{affine} monoids and their monoid rings. In Section
\ref{pyra} we describe a geometric reduction in Conjecture
\ref{conj} -- \emph{pyramidal descent}. It is a basic
combinatorial tool that puts the polyhedral/combinatorial
complexity of the nilpotence conjecture in the paradigm of the
existing algebraic techniques. An arithmetic analog of the main
result of \cite{Gu6} is developed in Section \ref{almost}. Section
\ref{homology} contains an interpretation of the arithmetic case
of the problem in terms of Hochschild homology. In Section
\ref{setup} we set up the notation and prove several lemmas to be
used in Section \ref{finite}. The latter section is devoted to the
solution to the aforementioned homological problem. In Section
\ref{essentially} we show how the arithmetic case implies the
general case.

\subsection{Conventions.}\label{conventions} The rings below,
unless specified otherwise, are assumed to be commutative and with
unit. Throughout the paper ``a sequence ${\bf
c}=(c_1,c_2,\ldots)$'' will \emph{always} mean a sequence of
natural numbers $\geq2$.

As usual, $\NN=\{1,2,\ldots\}$, $\ZZ_+=\{0,1,2,\ldots\}$, and
$\RR_+$ refers to the nonnegative real numbers.

\subsection{Acknowledgment.}\label{thanks}
The present exposition differs substantially from the first
version of the paper. I am extremely grateful to the referee for
the careful reading of the text and very insightful suggestions
that led me to many changes in the exposition, shortcuts (in
Section \ref{homology}), additions (in Section \ref{setup}), and
corrections of inaccuracies (in Section \ref{finite}).

\section{Polytopes, cones, and monoids}\label{preliminaries}

For the proofs of the general facts on polytopes and cones we use
below the reader is referred to \cite{O} and \cite{Z}. For a
systematic theory of the interplay between affine monoids and
polyhedral geometry see \cite{BrGu}.

\subsection{\emph{Polytopes}}\label{polytopes}
For a subset $W$ of an Euclidean space $\RR^n$ by $\conv(W)$ we
denote the convex hull of $W$, i.~e. the minimal convex subset of
$\RR^n$ containing $M$:
\begin{align*}
\conv(W)=\{a_1x_1+\cdots+a_kx_k\ :\ k\in\NN,\ x_1,\ldots,x_k\in
W,\ a_1,\ldots,a_k\in\RR,\\ 0\le a_1,\ldots,a_k\le1,\
a_1+\cdots+a_k=1\}.
\end{align*}

If we drop the condition  $0\le a_1,\ldots,a_k\le1$ what we get is
called the \emph{affine hull} of $W$ and denoted by $\aff(W)$.

The convex hull of a finite set is called a \emph{polytope}.
Polytopes are exactly the bounded sets which can be represented as
the intersections of finite systems of closed half-spaces. The
\emph{dimension} of a polytope $P$ is that of the affine hull of
$P$.

A polytope is a \emph{simplex} if the number of its vertices
equals its dimension $+1$.

If $P\subset\RR^n$ is a polytope and $\mathcal H\subset\RR^n$ is a
half-space, containing $P$, whose boundary is a hyperplane
$\partial\mathcal H$, then $P\cap\partial\mathcal H$ is called a
\emph{face} of $P$. The $0$-dimensional faces of $P$ are called
\emph{vertices}, $1$-dimensional faces are called \emph{edges},
and $(\dim P-1)$-dimensional (that is, codimension 1) faces are
called \emph{facets}.

For a polytope $P$ the union of its proper faces is called the
\emph{boundary} of $P$. It is denoted by $\partial P$. The subset
$P\setminus\partial P$ is called the \emph{relative interior} of
$P$ and is denoted by $\inte(P)$. Since we will only consider
relative interiors the adjective `relative' will be omitted.

The \emph{combinatorial type} of a polytope is the vertex-face
incidence matrix. If two polytopes have same vertex-face incidence
matrices up to a bijection between their vertex sets then we say
that the polytopes are \emph{of the same combinatorial type}. Two
polytopes of the same combinatorial type have the same number of
equidimensional faces. In particular, such polytopes have the same
dimension.

A polytope $P\subset\RR^n$ is called \emph{rational} if its
vertices belong to $\QQ^n$.

\subsection{\emph{Cones}}\label{conees}
A subset $C\subset\RR^n$ is called a \emph{convex cone} if $C$ is
a convex set and $ax\in C$ whenever $x\in C$ and $a\in\RR_+$.

A blanket assumption for the paper is that the cones considered
below are convex, positive dimensional, \emph{finite-polyhedral}
and \emph{pointed}. Here ``finite-polyhedral'' means ``the
intersection of a finite system of closed half-spaces whose
boundaries contain the origin $O\in\RR^n$'' and ``pointed'' means
``containing no affine line''.

A cone is called \emph{simplicial} if its edges are spanned by
linearly independent vectors.

As in the case of polytopes one can introduce the notion of a
face, the boundary, the interior, and the combinatorial type of a
cone.

For a cone its only vertex is the origin $O$ and the \emph{edges}
(the 1-dimensional faces) are the extremal rays containing the
origin $O$. A cone is the convex hull of its edges. In a sense,
the edges for a cone are the same as the vertices for a polytope.

We say that a cone is \emph{rational} if its edges contain
rational (equivalently, integral) points of $\RR^n$.

For a cone $C\subset\RR^n$ there exists a hyperplane
$H\subset\RR^n\setminus\{O\}$ such that the intersection
$\Phi(C)=H\cap C$ is a polytope. Moreover, $H$ can be chosen to be
rational. This follows, for instance, from the easily checked fact
that $C$ can be embedded into the positive orthant $\RR^n_+$. In
this situation $\RR_+\Phi(C):=\{ax\ :\ a\in\RR_+,\
x\in\Phi(C)\}=C$ and $\Phi(C)$ is a rational polytope whenever $C$
is a rational cone.

Of course, the polytope $\Phi(C)$ depends on the choice of $H$ and
is, therefore, only defined up to a projective transformation. But
the combinatorial types of $\Phi(C)$ and $C$ determine each other.
In particular, $C$ is simplicial if and only if $\Phi(C)$ is a
simplex.

When we consider the $\Phi$-polytopes simultaneously for several
cones contained in a single cone it is always assumed that these
polytopes live in the \emph{same} rational hyperplane $H$ that has
been chosen for the ambient cone.

\subsection{\emph{Monoids}.}\label{Monoids}
All monoids $M$ considered in this paper are assumed to be
commutative, cancellative and torsion free. In other words, we
assume that the natural homomorphisms
$M\to\gp(M)\to\QQ\otimes\gp(M)$ are injective, where $\gp(M)$
refers to the universal group associated to $M$ (the {\it
Grothendieck group}, or the {\it group of differences} of $M$).
Equivalently, our monoids can be thought of as additive submonoids
of real spaces ($M\to\RR\otimes\gp(M)$). This enables us to use
polyhedral geometry in their study. Our conditions on a monoid $M$
are equivalent to the condition that $\k[M]$ is a domain for some
(equivalently, arbitrary) field $\k$.

We put $\rank M=\rank\gp(M)=\dim_{\QQ}\QQ\otimes\gp(M)$.

A finitely generated monoid will be called an \emph{affine
monoid}.

When we treat monoids separately, i.~e. not within their monoid
rings, the monoid operation is written additively. In monoid rings
we switch to the multiplicative notation.

For a monoid $M$ its group of units (i.~e. the maximal subgroup of
$M$) is denoted by $U(M)$. Monoids for which $U(M)=0$ are called
\emph{positive monoids}. Any affine positive monoid embeds into a
free monoid $\ZZ_+^r$ for some $r\in\NN$ (\cite[\S2.1]{BrGu}).
Therefore, for any ring $R$ and any affine positive monoid $M$
there is a grading
$$
R[M]=R\oplus R_1\oplus\cdots
$$
making the elements of $M$ homogeneous.

For a nonzero affine positive monoid $M\subset\ZZ^n$ we put
$$
M_*=\bigl(\inte(\conv(M))\cap M\bigr)\cup\{O\},
$$
$\conv(M)$ being considered in $\RR\otimes\gp(M)$. In algebraic
terms, $M_*=\{m\in M\ :\ \ZZ m+M=\gp(M)\}\cup\{O\}$ (see
\cite[Section 2.1, \emph{Inversion of monomials}]{BrGu}). The
monoid $M_*$ is never affine, unless $\rank M=1$ in which case
$M_*=M$.

An affine positive monoid $M$ defines the cone $\RR_+\otimes
M\subset\RR\otimes\gp(M)$ (i.~e. $\RR_+\otimes
M=\conv(\I(M\to\RR\otimes\gp(M))$). We put
$\Phi(M):=\Phi(\RR_+\otimes M)$. We say that $M$ is
\emph{simplicial} if $\Phi(M)$ is a simplex.

\subsection{\emph{Normal monoids}}\label{Normality}
A monoid $M$ is called \emph{normal} if the submonoid $\bar
M:=\{x\in\gp(M)\ :\ nx\in M\ \text{for some natural}\
n\}\subset\gp(M)$ coincides with $M$. In general, the monoid $\bar
M$ may be larger than $M$. We call it the \emph{integral closure}
of $M$. It is a well known fact that a monoid $M$ is normal if and
only if the monoid algebra $\k[M]$ is integrally closed domain for
some (equivalently, arbitrary) field $\k$.

Consider an affine positive monoid $M$. Upon fixing an embedding
$M\to\RR\otimes\gp(M)\cong\RR^n$, $n=\rank(M)$ and making the
identification $\gp(M)=\ZZ^n$ we can view the monoid $M$ as a set
of integral points in $\RR^n$. The convex hull $C=\conv(M)$ is
then an $n$-dimensional rational cone in $\RR^n$ and the normality
condition translates into the equality $M=C\cap\ZZ^n$.

Conversely, for a cone $C\subset\RR^n$ the submonoid
$C\cap\ZZ^n\subset\ZZ^n$ is a normal affine positive monoid. The
finite generation part of this claim is the classical \emph{Gordan
Lemma} \cite{Gor}.

Next is an extract from \cite{Gu6}.

\begin {proposition}\label{norint} Let $R$ be a regular ring and
$p$ be a natural number. Then the following are true:

(a) $\NN$ acts nilpotently on $K_p(R[M])$ for every monoid $M$ if
this is so for affine normal positive monoids.

(b) Assume $\NN$ acts nilpotently on $K_p(R[M])$ for every affine
normal positive monoid $M$ of rank $<n$. If $N$ is a rank $n$
affine positive monoid then for every sequence ${\bf
c}=(c_1,c_2,\ldots)$ and every element $x\in K_p(R[N])$ there
exists $j_0$ such that $(c_1\cdots
c_j)_*(x)\in\I\bigl(K_p(R[N_*])\to K_p(R[N])\bigr)$ for $j>j_0$.

(c) Conjecture \ref{conj} is true for all affine simplicial
monoids.
\end{proposition}

In fact, (a) and (b) together constitute a reformulation of
\cite[Lemma 3.4]{Gu6} and the observations on $K$-theoretical
excision in \cite[\S3]{Gu6}; (c) is \cite[Theorem 6.4]{Gu6}.

\section{Sufficiency of pyramidal descent}\label{pyra}

In this section we make a combinatorial reduction in Conjecture
\ref{conj}, called \emph{pyramidal descent}.

\subsection{\emph{Pyramidal extensions}}\label{Pyramidal}
A polytope $P\subset\RR^n$ is called a \emph{pyramid} if it is a
convex hull of one of its facets $F\subset P$ and a vertex $v\in
P$, not in the affine hull of $F$. In this situation $F$ is a
\emph{base} and $v$ is an \emph{apex} of $P$, and we write
$P=\pyr(v,F)$. For instance, an arbitrary simplex is a pyramid
such that every facet is a base and every vertex is an apex.

The \emph{complexity} of a $d$-dimensional polytope
$P\subset\RR^n$ ($d\leq n$) is by definition the number
$\cc(P)=d-i$, where $i$ is the maximal nonnegative integer
satisfying the condition: there exists a sequence $P_0\subset
P_1\subset\cdots\subset P_i=P$ such that $P_j$ is a pyramid over
$P_{j-1}$ for each $1\leq j\leq i$.

Thus the complexity of a polytope is measured by the number of
steps needed to get to the polytope by successively taking
pyramids over some initial polytope: the more steps we need the
simpler the polytope is. It is immediately observed that: (i) the
complexity is an invariant of the combinatorial type, (ii) a
polytope $P$ is not a pyramid if and only if $\cc(P)=\dim P$,
(iii) simplices are exactly the polytopes of complexity $0$, and,
(iv) we always have the equality $\cc(\pyr(v,P))=\cc(P)$.

For a cone $C\subset\RR^n$ its combinatorial complexity $\cc(C)$
is defined to be $\cc(\Phi(C))$, and for a positive affine monoid
$M$ its combinatorial complexity $\cc(M)$ is defined to be
$\cc(\Phi(M))$.

Consider two polytopes $P\subset Q$, $P\not=Q$. Assume $P$ is
obtained from $Q$ by cutting off a pyramid at a vertex $v\in Q$.
In other words, $Q=P\cup P'$, $\dim P=\dim P'=\dim Q$ and
$P'=\pyr(v,P\cap P')$. In this situation we say that $P\subset Q$
is a \emph{pyramidal extension}. Observe, that $P\subset Q$ is a
pyramidal extension of polytopes then these polytopes are
automatically positive dimensional.

An extension of monoids $M\subset N$ is called \emph{pyramidal} if
$M$ and $N$ are nonzero affine positive normal monoids,
$\gp(M)=\gp(N)$, and $\Phi(M)\subset\Phi(N)$ is a pyramidal
extension of polytopes.

The following lemma is a key combinatorial fact. Let
$P\subset\RR^n$ be a polytope. Call a sequence of polytopes
$P=P_0,P_1,P_2,\ldots$ \emph{admissible} if the following
conditions hold for all indices $k$:
\begin{itemize}
\item[(1)] $P_k\subset P$, \item[(2)] either $P_{k+1}\subset P_k$
is a pyramidal extension or $P_k\subset P_{k+1}$.
\end{itemize}

(Observe, $\dim P_k=\dim P_0$ for all $P_k$ in an admissible
sequence $P_0,\ldots,P_k,\ldots$)

\begin{lemma}[\cite{Gu1}]\label{admis}
For any open subset $U$ of a polytope $P$ there exits an
admissible sequence of polytopes $P=P_0,P_1,P_2,\ldots$ such that
$P_j\subset U$ for all sufficiently big indices $j$. If $P$ is
rational then the polytopes $P_j$ can be chosen to be rational.
\end{lemma}

Let $M\subset N$ be a pyramidal extension of monoids. It will be
called \emph{an extension of complexity} $c$ if
$\cc\bigl(\overline{\Phi(N)\setminus\Phi(M)}\bigr)=c$, where
$\overline Y$ refers to the Euclidean closure of a subset
$Y\subset\RR^n$. In this situation we will write $\cc(M\subset
N)=c$.

\subsection{\emph{Pyramidal descent}}\label{Descent}
Assume we are given a regular ring  $R$, a natural number $p$, and
a pyramidal extension of monoids $L\subset M$. We say that
\emph{pyramidal descent} holds for these objects if for every
element $x\in K_p(R[M])$ and every sequence ${\bf
c}=(c_1,c_2,\ldots)$ there exists an index $j_0\in\NN$ such that
$(c_1\cdots c_j)_*(x)\in\I{\big(}K_p(R[L])\to K_p(R[M]){\big)}$
for all $j>j_0$. We say that \emph{pyramidal descent of type $c$
holds} (for $R$ and $p$) if pyramidal descent holds for the
$K_p$-groups of the monoid $R$-algebras, corresponding to all
pyramidal extensions of monoids $L\subset M$ where $\cc(L\subset
M)=c$. If we require the latter condition only for rank $r$
monoids we say that pyramidal descent of type $c$ holds \emph{for
monoids of rank $r$}. Pyramidal descent, without referring to the
complexity and the rank, just means the corresponding condition
for all pyramidal extensions (for $R$ and $p$ fixed). Also, if $R$
and $p$ are clear from the context we will usually skip them.

Next is a refined version of Lemma 5.2 in \cite{Gu6}.

\begin{proposition}\label{sufficient}
Let $R$ and $p$ be as above. Assume $N$ is an affine normal
positive monoid of positive complexity. Then Conjecture \ref{conj}
is true for the monoid ring $R[N]$ if pyramidal descent of type
$<\cc(N)$ holds for monoids of rank $\rank(N)$.
\end{proposition}

\begin{proof}
Let $P_0\subset P_1\subset\cdots\subset P_i=\Phi(N)$ be a sequence
of polytopes where $i=\rank(N)-1-\cc(N)$ ($=\dim\Phi(N)-\cc(N)$)
and $P_j=\pyr(v_j,P_{j-1})$ for each $j\in[1,i]$.

Fix a rational simplex $\Delta\subset P_0$, $\dim\Delta=\dim P_0$.
By Lemma \ref{admis} there is an admissible sequence
$P_0=Q_0,Q_1,Q_2,\ldots$ of rational polytopes such that
$Q_t\subset\Delta$ for all $t\gg0$. Then the sequence of polytopes
$\tilde Q_t=\conv(v_1,\ldots,v_i,Q_t)$ is an admissible sequence
of rational polytopes such that $\tilde Q_0=\Phi(N)$ and $\tilde
Q_t$ are contained in the simplex
$\tilde\Delta=\conv(v_1,\ldots,v_i,\Delta)$ for $t\gg0$. (We
assume $\tilde Q_t=Q_t$ and $\tilde\Delta=\Delta$ when $i=0$)
Moreover, if $Q_{t+1}\subset Q_t$ is a pyramidal extension then we
have $\cc\bigl(\overline{\tilde Q_{t+1}\setminus\tilde
Q_t}\bigr)\leq\cc(N)-1$. By Gordan's Lemma the monoids
$\RR_+\tilde Q_t\cap\gp(N)$ are all affine. Obviously, they are
also normal and positive.

Fix a sequence ${\bf c}=(c_1,c_2,\ldots)$ and an element $x\in
K_p(R[N])$. Assume pyramidal descent of type $\leq\cc(N)-1$ holds
for rank $r$ monoids. Then there exists a sequence of elements
$$
x_t\in K_p(R[\RR_+\tilde Q_t\cap\gp(N)]),\quad t\in\ZZ_+
$$
and that of nonnegative integers $0=k_0\leq k_1\leq k_2\leq\cdots$
such that the following conditions hold:
\begin{itemize}
\item $x_0=x$, \item if $\tilde Q_{t+1}\subset\tilde Q_t$ is a
pyramidal extension then $k_t<k_{t+1}$ and $x_{t+1}$ maps to
$(c_{k_t+1}c_{k_t+2}\cdots c_{k_{t+1}})_*(x_t)$ under the map
$$
K_p(R[\RR_+\tilde Q_{t+1}\cap\gp(N)])\to K_p(R[\RR_+\tilde
Q_t\cap\gp(N)]),
$$
\item if $\tilde Q_t\subset\tilde Q_{t+1}$ then $k_t=k_{t+1}$ and
$x_{t+1}$ is the image of $x_t$ under the map $K_p(R[\RR_+\tilde
Q_t\cap\gp(N)])\to K_p(R[\RR_+\tilde Q_{t+1}\cap\gp(N)])$.
\end{itemize}

In this situation we have
$$
(c_1\cdots
c_j)_*(x)\in\I\bigl(K_p(R[\RR_+\tilde\Delta\cap\gp(N)])\to
K_p(R[N])\bigr)
$$
for all $j\gg0$ and we are done by Proposition \ref{norint}(c).
\end{proof}

\begin{remark}\label{Admissible}
(a) The main $K$-theoretical difficulty we face here is to achieve
the pyramidal descent. However, we also need to allow the
containment $Q_t\subset Q_{t+1}$ in our definition of admissible
sequences, trivial from the $K$-theoretical point of view. Only
resorting to `removing' vertices from polytopes would lead to
increasingly complex combinatorial types, preventing us from
contracting polytopes into arbitrarily small neighborhoods.

(b) It follows from Proposition \ref{sufficient} that for rank 3
monoids we only need to achieve pyramidal descent of type 0.
\end{remark}

By induction first on complexity and then on rank, Propositions
\ref{norint}(a,c) and \ref{sufficient} imply the following
proposition, in the formulation of which we assume that a regular
ring $R$ and a $K$-theoretical index $p$ are fixed and `monoid'
refers to an affine, positive, normal monoid.

\begin{proposition}\label{induction}
To prove Conjecture \ref{conj} it is enough to show that the two
things:
\begin{itemize}
\item
\noindent Conjecture \ref{conj} for monoids of rank $<n$, and
\item
Conjecture \ref{conj} for monoids of rank $n$ and complexity $<c$
\end{itemize}
together imply pyramidal descent of type $<c$ for monoids of rank
$n$.
\end{proposition}

\begin{remark}\label{rankcomp}
The implication, mentioned in Proposition \ref{induction}, has
`almost been shown' in \cite{Gu6}, see the next section. The
reason we introduced the notion of complexity for pyramidal
extensions is that it allows a more subtle passage from rank $<n$
monoids to rank $n$ monoids. Without such an approach the
inductive step on the monoid rank would raise the transcendence
degree of the coefficient field and, by the same token, we would
run into the class of non-arithmetic fields.
\end{remark}

\section{Almost pyramidal descent}\label{almost}

\subsection{\emph{Finitely many exceptional
monomials}}\label{newdescent}
In this section we fix a field $\k$ of characteristic 0 and a
natural number $p$.

We now develop a version of the main result of \cite{Gu6}, adapted
to the arithmetic situation.

Consider the following construction. Let $M\subset N$ be a
pyramidal extension of rank $n$ monoids. We make the
identification $\gp(N)=\ZZ^n$ and, hence,
$\RR\otimes\gp(N)=\RR^n$. Assume we are given a cone
$D\subset\RR_+M_*$ ($=\conv(M_*)$), a rational point
$v\in\Phi(N)\setminus\Phi(M)$, and a nonzero element
$t\in(\RR_+v)\cap\ZZ^n\cong\ZZ_+$. To these objects we associate
the following noncommutative ring
$$
\Lambda=\Lambda(D,t)=
\begin{pmatrix}
\k[L]&\k[L]\cap t^{-1}\k[L]\\
\k[L]+t\k[L]& \k[L]
\end{pmatrix}
\subset M_{2\times 2}(\k[N]),
$$
where $L=D\cap\ZZ^n$ and $M_{2\times 2}(\k[N])$ is the ring of
$2\times2$-matrices over $\k[N]$.

\begin{theorem}[Almost Pyramidal Descent]\label{maing6}
Assume Conjecture \ref{conj} holds for the monoid $\k$-algebras of
all affine normal positive monoids of rank $<n$ and assume it
holds for all affine normal positive monoids of rank $n$ and
complexity $<\cc(N)$. Assume further $\cc(M\subset N)<\cc(N)$.
Then for every sequence ${\bf c}=(c_1,c_2,\ldots)$ and every
element $x\in K_p(\k[N])$ there exist a cone $D\subset\RR_+M_*$, a
rational point $v\in\Phi(N)\setminus\Phi(M)$, and a nonzero
element $t\in(\RR_+v)\cap\ZZ^n$, such that:
$$
(c_1\cdots c_j)_*(x)\in\I\bigl(K_p(\Lambda)\to K_p(M_{2\times
2}(\k[N]))=K_p(\k[N])\bigr)
$$
for all $j\gg0$.
\end{theorem}

\begin{remark}\label{comments}
(a) We call this theorem `almost pyramidal descent' (over $\k$)
because of the following. The ring $\Lambda$ is generated as a
$\k$-algebra by monomials, i.~e. by $2\times2$-matrices whose
entries are all zero but one and the distinguished entry belongs
to $N$. The set of monomials of $\Lambda$ not in the matrix
subring $M_{2\times2}(\k[M])\subset M_{2\times2}(\k[N])$ is
\emph{finite} -- this is so because $D\subset\RR_+M_*$. The goal
is to show that these \emph{exceptional monomials} do not
contribute to the $K$-theory of $\k[N]$ in the following sense:
$$
(c_1\cdots c_j)_*(x)\in\I{\big(}K_p(\Lambda\cap
M_{2\times2}(\k[M]))\to K_p(M_{2\times
2}(\k[N]))=K_p(\k[N]){\big)},\ j\gg0.
$$
This would imply pyramidal descent for the extension $\k[M]\subset
\k[N]$. In other words, we only need to show that the exceptional
monomials can always be eliminated.

(b) Actually,  \cite[Theorem 9.3]{Gu6}, on which Theorem
\ref{maing6} is based, is proved for fields of any characteristic.
It even says that there is \emph{only one} monomial to be
eliminated. But what is really relevant to the issue is the
finiteness of the set of exceptional monomials.
\end{remark}

A word on one notation. For any (not necessarily affine) submonoid
$H\subset N$ we have the $\k$-algebra augmentation $\k[H]\to \k$,
$H\setminus\{1\}\to0$, whose right inverse is the identity
embedding $\k\to \k[H]$. In particular, the homomorphism
$K_p(\k)\to K_p(\k[H])$ is a split monomorphism and, therefore,
$K_p(\k)$ can be thought of as a natural direct summand of
$K_p(\k[H])$. For an element $\tau\in K_p(\k[H])$ its image in
$K_p(\k)$ will be denoted by $\tau(0)$.

Consider two submonoids $H_1\subset H_2\subset N$ and an element
$\tau_1\in K_p(\k[H_1])$. Denote by $\tau_2$ the image of $\tau_1$
in $K_p(\k[H_2])$. Then
\begin{equation}\label{augmentation}
\tau_1(0)=0\ \iff\ \tau_2(0)=0.
\end{equation}

\subsection{\emph{Previous result}}\label{previous}
Before embarking on the proof of Theorem \ref{maing6} we recall
\cite[Theorem 9.3]{Gu6} and the part of its proof to be adapted to
the new situation.

Fix a sequence ${\bf c}=(c_1,c_2,\ldots)$ and an element $z\in
K_p(\k[N_*])/K_p(\k)$. \cite[Theorem 9.3]{Gu6} asserts that there
exists a cone $D\subset\RR_+M_*$, a rational point
$v\in\inte(\Phi(N))\setminus\Phi(M)$, a nonzero element
$t\in(\RR_+v)\cap\ZZ^n$, and a normal submonoid $\tilde L\subset
D\cap\ZZ^n$ with the property $\RR_+\tilde L=D$, such that
$$
(c_1\cdots c_j)_*(z)\in\I{\big(}K_p(\Lambda')\to K_p(M_{2\times
2}(\k[N_*]))=K_p(\k[N_*]){\big)},\quad j\gg0,
$$
where
$$
\tilde\Lambda=
\begin{pmatrix}
\k[\tilde L]&\k[\tilde L]\cap t^{-1}\k[\tilde L]\\
\k[\tilde L]+t\k[\tilde L]& \k[\tilde L]
\end{pmatrix}
\subset M_{2\times 2}(\k[N_*]).
$$
Here $\tilde L$ and $\tilde\Lambda$ are the same as $L(\Gamma)$
and $\Lambda_{t,\Gamma,L}$ in \cite[Theorem 9.3]{Gu6} ($\Gamma$
being a subpolytope in the interior of $\Phi(M)$).

In fact, the only difference between the assumptions in Theorem
\ref{maing6} and the assumptions in \cite[Theorem 9.3]{Gu6} is
that in \cite{Gu6} we require that Conjecture \ref{conj} holds for
the monoid algebras of all rank $<n$ monoids without nontrivial
units and for arbitrary coefficient fields. On the other hand
there is no mention in \cite[Theorem 9.3]{Gu6} about the
complexity of the pyramidal extension $M\subset N$ and the proof
of \cite[Theorem 9.3]{Gu6} only refers to affine positive normal
monoids.\footnote{In addition, the monoid $\tilde L$ is of very
special type -- s.~c. \emph{polarized monoid}}

\

It is essential to notice that the proof in \cite{Gu6} only uses
the fields $\k$ and $\k(T)$, $T$ a variable. Assume we have shown
how to circumvent the use of the field $\k(T)$ in the situation of
Theorem \ref{maing6}. By Proposition \ref{norint}(b) there is no
loss of generality in assuming that $x\in\I\bigl(K_p(\k[N_*])\to
K_p(\k[N])\bigr)$. Replacing $x$ by $x-x(0)$ we can further assume
that $x$ is the image of an element $z\in K_p(\k[N_*])/K_p(\k)$.
Therefore, Theorem \ref{maing6} follows by extending $\tilde L$ to
the monoid $L=D\cap\ZZ^n$ and extending $\tilde\Lambda$ to the
ring $\Lambda=\Lambda(D,t)\subset M_{2\times 2}(\k[N_*])$.

\

Next we describe where and how the field $\k(T)$ is used in the
proof of \cite[Theorem 9.3]{Gu6}. It is used in Step 2 in that
proof, where we have:
\begin{itemize}
\item[(i)] an arbitrary element $z\in K_p(\k[N_*])/K_p(\k)$,
\item[(ii)] a cone $D_z\subset\RR_+M_*$, \item[(iii)] a rational
point $v_z\in\inte(\Phi(N))\setminus\Phi(M)$, \item[(iv)] a
nonzero element $t_z\in\RR_+v_z\cap\ZZ^n$, \item[(v)] a submonoid
$L_z\subset D_z\cap\ZZ^n$,
\end{itemize}
such that the following conditions hold:
\begin{itemize}
\item[(vi)] $\RR_+L_z=D_z$, \item[(vii)] the submonoid
$\ZZ_+t_z+L_z\subset N_*$ is normal, \item[(viii)]
$z\in\I\bigl(K_p(\k[\ZZ_+t_z+L_z])\to K_p(\k[N_*])\bigr)$.
\end{itemize}

Moreover, for any real number $\epsilon>0$ (and fixed $z$) the
objects $D_z$, $v_z$, $t_z$ and $L_z$ can be chosen in such a way
that:
\begin{itemize}
\item[(ix)]
$\partial\Phi(\ZZ_+t_z+L_z)=\partial\conv(v_z,\Phi(D_z))$ is
within the $\epsilon$-neighborhood of $\partial\Phi(N)$.
\end{itemize}
We explain the latter condition in detail. The element $t_z$, the
polytope $\Phi(D_z)$, and the monoid $L_z$ are respectively the
pole $t$, the polytope $\Gamma_1$, and the monoid $N_1(\Gamma_1)$
in the proof of \cite[Theorem 9.3]{Gu6}. In particular,
$\ZZ_+t_z+L_z=N_1$. It follows from the approximation principle
(Theorem 6.2) in \cite[Section 6]{Gu6} that the pole $t$ and the
subpolytope $\Gamma_1\subset\inte(\Phi(M))$ can be chosen in such
a way that $v=\Phi(t)$ and $\Gamma_1$ approximate respectively the
vertex of $\Phi(N)$, not in $\Phi(M)$, and the open set
$\inte(\Phi(M))$ with arbitrary prespecified precision.
Informally, the better approximation we take the better the proof
of [Gu6,Theorem 9.3] works.

\

Since $\ZZ_+t_z+L_z$ is normal, its localization at $t_z$ (i.~e.
the submonoid $\ZZ t_z+L_z\subset\gp(N)$) is isomorphic to a
monoid of the type $\ZZ\oplus L'$ for a certain affine normal
positive monoid $L'$ of rank $n-1$, a general observation on
inversion of extremal elements in a normal positive monoid (see,
for instance, \cite[\S2.1]{BrGu}). Under this isomorphism $t_z$
maps to $(1,0)$. Furthermore, we can fix an isomorphism $\ZZ\oplus
L'\cong\ZZ t_z+L_z$ in such a way that the image of $\ZZ_+\oplus
L'$ contains $\ZZ_+t_z+L_z$, see \cite[Corollary 2.5]{Gu6}. We
will simply identify $\ZZ_+\oplus L'$ with its image in $\ZZ
t_z+L_z$.

Because of the condition (ix) we can assume that
$$
z\in\I\bigl(K_p(\k[\ZZ_+t_z+L_z])\to K_p(\k[N_*])\bigr).
$$
\noindent Let $z_1$ be a preimage of $z$ and let $z_2$ be the
image of $z_1$ in $K_p(\k[\ZZ_+t_z+L'])$. The field $\k(T)$, being
identified with the fraction field of $\k[\ZZ_+t_z]$, is used in
\cite{Gu6} to show that:
\begin{itemize}
\item[(x)] \ \ $(c_1\cdots c_j)_\bullet(z_2)=0$ for $j\gg0$,
\end{itemize}
\noindent where for a natural number $c$ by $c_\bullet$ we denote
the endomorphism of the group $K_p(\k[\ZZ_+t_z+L'])$ induced by
the $\k[\ZZ_+t_z]$-algebra endomorphism of $\k[\ZZ_+t_z+L']$,
$l\mapsto l^c$, $l\in L'$. Notice, up to the notation change
$T\mapsto t$, $\bf c\to{\mathfrak c}'$, $z_2\mapsto x_1$, the
condition (x) is the same as the inclusion
$x_1\in\Ker\left(K_p(\k[N_1])\to K_p\left(\k[N_0^{{\bf
c}'}[t]\right)\right)$ at the end of page 310 in \cite{Gu6}.

In particular, Theorem \ref{maing6} is proved once we know how to
achieve the condition (x) above, working \emph{exclusively} over
the field $\k$. Since the choice of the preimage $z_1$ is
irrelevant in the proof of \cite[Theorem 9.3]{Gu6}, this will be
achieved by a more careful choice with use of the condition on
complexities.

\subsection{\emph{Proof of Theorem \ref{maing6}}}\label{proof4.1}
We keep the notation in Section \ref{previous}. Let $w$ denote the
vertex of $\Phi(N)$ not in $\Phi(M)$. We fix arbitrarily a
rational point
$$
\xi\in\inte\big(\Phi(M)\cap\bigl(\overline{\Phi(N)
\setminus\Phi(M)}\bigr)\bigr).
$$
In other words, $\xi$ is an interior rational point of the facet
of $\Phi(M)$ visible from $w$. For a \emph{rational} number
$0<\lambda<1$ we let $\Phi(N)^\lambda$ denote the homothetic image
of $\Phi(N)$, centered at $\xi$ and with factor $\lambda$. Put
$N^\lambda=\RR_+\Phi(N)^\lambda\cap\ZZ^n$. Then $N^\lambda$ is an
affine normal positive submonoid of $N_*$. (The finite generation
of $N^\lambda$ follows from Gordan's lemma, see Section
\ref{Normality}.)

By first choosing a rational number $0<\lambda<1$ close enough to
1 and then choosing correspondingly a real number $\epsilon>0$
small enough we can achieve that:
\begin{itemize}
\item $z\in\I\bigl(K_p(\k[N^\lambda])\to K_p(\k[N_*])\bigr)$,
\item $\Phi(N)^\lambda\subset\Phi(\ZZ_+t_z+L_z)$.
\end{itemize}
Let $w^\lambda$ denote the image of $w$ under our homothetic
transformation. The intersection of the cone
$w^\lambda+\RR_+\left(\Phi(N)^\lambda-w^\lambda\right)
\subset\aff(\Phi(N))$
and the pyramid $\Phi(\ZZ_+t_z+L')$ is a rational polytope. Denote
it by $P$ and consider the monoid $N(P)=\RR_+P\cap\ZZ^n$. We have:
\begin{itemize}
\item if $\epsilon>0$ is small enough then $P$ is a rational
pyramid with apex at $w^\lambda$, having the same combinatorial
type as $\overline{\Phi(N)\setminus\Phi(M)}$, \item
$N^\lambda\subset N(P)\subset\ZZ_+t_z+L'$.
\end{itemize}
We fix $\lambda$ and $\epsilon$ in such a way that all these
conditions are satisfied.

Let $z_0\in K_p(\k[N^\lambda])$ map to $z$ under the homomorphism
$K_p(\k[N^\lambda])\to K_p(\k[N_*])$. Since $z(0)=0$ we have
$z_0(0)=0$ by (\ref{augmentation}).

Let $z_1\in K_p(\k[\ZZ_+t_z+L_z])$, $z_2\in K_p(\k[\ZZ_+t_z+L'])$,
and $z_3\in K_p(\k[N(P)])$ be the corresponding images of $z_0$.
Then $z_1(0)=0$, $z_2(0)=0$ and $z_3(0)=0$ by
(\ref{augmentation}).

Since $\cc(N(P))=\cc(M\subset N)$ and $z_3(0)=0$, by the
assumption in Theorem \ref{maing6} we have $(c_1\cdots
c_j)_*(z_3)=0$ for $j\gg0$. (As usual, for a natural number $c$
the endomorphism of $K_p(\k[N(P)])$, induced by the monoid
endomorphism $N(P)\to N(P)$, $n\mapsto n^c$, is denoted by $c_*$.)

Since $z_2$ is the image of $z_3$ it follows that $(c_1\cdots
c_j)_*(z_2)=0$ in $K_p(\k[\ZZ_+t_z+L'])$ whenever $j\gg0$. Here we
use the same notation $c_*$ ($c\in\NN$) for the endomorphism of
$K_p(\k[\ZZ_+t_z+L'])$, induced by $l\mapsto l^c$,
$l\in\ZZ_+t_z+L'$.

Therefore, to prove (x) for our element $z_2$ it only remains to
show that
\begin{equation}\label{zeta}
c\in\NN,\ \zeta\in K_p(\k[\ZZ_+t_z+L'])/K_p(\k),\ c_*(\zeta)=0\
\Longrightarrow\ c_\bullet(\zeta)=0.
\end{equation}
For a natural number $c$ we have the equality $c_*=c^\bullet
c_\bullet$ of the endomorphisms of $K_p(\k[\ZZ_+t_z+L'])$, where
$c_\bullet$ is as in (x) in Section \ref{previous} and $c^\bullet$
is induced by the $\k[L']$-algebra endomorphism of
$\k[\ZZ_+t_z+L']$, $t_z\mapsto t_z^c$. This algebra endomorphism
makes $\k[\ZZ_+t_z+L']$ a $\rank c$ free module over itself.
Therefore, by the transfer map for $K_p$ we get $c\cdot\zeta=0$
for every $c\in\NN$ and $\zeta\in K_p(\k[\ZZ_+t_z+L'])/K_p(\k)$
such that $c_*(\zeta)=0$.

Now the implication (\ref{zeta}) follows from a graded structure
$\k[\ZZ_+t_z+L']=\k\oplus R_1\oplus\cdots$ (see Section
\ref{Monoids}) and the general fact which we borrow from the next
section (see Proposition \ref{totaro}(a)): for a not necessarily
commutative graded $\k$-algebra $A=A_0\oplus A_1\oplus\cdots$ and
the ideal $A^+=0\oplus A_1\oplus A_2\oplus\cdots\subset A$ the
relative group $K_p(A,A^+)$ is a $\k$-vector space. \qed

Notice, we use the condition $\chara \k=0$ only at the very end of
the proof above.

By Proposition \ref{induction} and Theorem \ref{maing6} we get

\begin{proposition}\label{descent}
In order to prove Conjecture \ref{conj} for monoid $\k$-algebras
it is enough to show that for any pyramidal extension of monoids
``almost pyramidal descent'' implies the actual pyramidal descent.
\end{proposition}

\section{Homological interpretation}\label{homology}

\subsection{\emph{Linear structures}}\label{Linear}
First we summarize general facts on the action of the big Witt
vectors on nil-$K$-groups -- the Bloch-Stienstra operations on
Nil-$K$-theory, generalized by Weibel to the graded noncommutative
situation \cite{Bl,St,W}.

In this subsection we assume $\k$ is a field of characteristic $0$
and $s\in\ZZ_+$.

For a not necessarily commutative graded $\k$-algebra $A=A_0\oplus
A_1\oplus\cdots$ ($\k$ is in the center of $A$) its ideal
generated by the elements of positive degree will be denoted by
$A^+$. Let $A_{[s]}=A_0\oplus0\oplus\cdots\oplus0\oplus A_s\oplus
A_{s+1}\oplus\cdots$. (In particular, $A=A_{[0]}=A_{[1]}$.)

In the special case $A=A_0[T]=A_0+A_0T+A_0T^2+\cdots$ ($T$ a
variable) the non-graded endomorphism of $A$, which is the
identity on $A_0$ and maps $T$ to $T^s$, induces a group
endomorphism $NK_p(A_0)\to NK_p(A_0)$ which is denoted by $\V_s$
and called the \emph{$s$th Verschiebung}.

The topological ring $\W(\k)$ of big Witt vectors over $\k$ is
isomorphic (through the \emph{ghost map}) to $\Pi_1^\infty k$. In
particular, $\k\subset\W(\k)$.

The interested reader is referred to \cite[Section 10]{Gu6} and
\cite{Gu8} for a brief survey of the following extract from
\cite{Bl,St,W}.

For a graded $\k$-algebra $A=A_0\oplus A_1\oplus\cdots$ we have
the multiplicative action
$$
\k\times A\to A,\quad (r,a_0+a_1+a_2+\cdots)\mapsto
a_0+ra_1+r^2a_2+\cdots
$$
which gives rise to an action
$$
\k\times K_p(A,A^+)\to K_p(A,A^+),\quad(r,w)\mapsto rw.
$$

\begin{proposition}\label{totaro} For each nonnegative integer $p$
the following hold:
\begin{itemize}
\item[(a)] The group $K_p(A,A^+)$ is a continuous $\W(\k)$-module;
as such it decomposes into eigenspaces for the homothety operators
$[r]$, $r\in\k$, the $i$th eigenspace being defined by the
condition $[r]w=r^iw$. This decomposition
$K_p(A,A^+)=\bigoplus_{j=1}^{\infty} V_j$ into $\k$-vector spaces
is functorial with respect to graded $\k$-algebra homomorphisms.
\item[(b)] For a decomposition as in (a) one has the inclusions
\begin{align*}
\I{\big(}&K_p(A_{[s]},A_{[s]}^+)\to K_p(A,A^+{\big)}\subset\\
&0\oplus\cdots\oplus0\oplus V_s\oplus V_{s+1}\oplus\cdots,\quad
s\in\NN.
\end{align*}
\item[(c)] The $s$th \emph{Verschiebung} $\V_s$ is a $\k$-linear
endomorphism of $NK_p(A_0)$, $s\in\NN$.
\end{itemize}
\end{proposition}

Let $M$ be an arbitrary affine positive monoid. Let
$\k[M]=B=\k\oplus B_1\oplus B_2\oplus\cdots$ be a grading, making
the elements of $M$ homogeneous (Section \ref{Monoids}).

\begin{lemma}\label{referee}
For any natural number $c$ the endomorphism $c_*:K_p(B,B^+)\to
K_p(B,B^+)$ is $\k$-linear.
\end{lemma}

\begin{proof}
Consider the commutative square of $\k$-algebras
\begin{equation}\label{linear}
\xymatrix{B[T]\ar[r]^{T\mapsto T^c}&B[T]\\
B\ar[r]_\alpha\ar[u]^\beta&B \ar[u]_\gamma}
\end{equation}
where $\alpha(m)=m^c$, $\beta(m)=T^{\deg(m)}m^c$,
$\gamma(m)=T^{\deg(m)}m$, $m\in M$. Thinking of $B[T]$ as a graded
ring with respect to the powers of $T$, the vertical maps in
(\ref{linear}) become graded homomorphisms. We arrive at the
commutative diagram of groups
$$
\xymatrix{NK_p(B)\ar[r]^{\V_c}&NK_p(B)\\
K_p(B,B^+)(\k)\ar[r]_{c_*}\ar[u]^{\beta_*}&K_p(B,B^+)
\ar[u]_{\gamma_*}}
$$
whose vertical maps are $\k$-linear by functoriality and the upper
horizontal map is $\k$-linear by Proposition \ref{totaro}(c). The
right vertical homomorphism is actually a monomorphism because the
map $T^{\deg(-)}(-)$ in diagram (\ref{linear}) splits the
(non-graded) $B$-algebra epimorphism $B[T]\to B$, $T\mapsto1$.
These conditions altogether imply the desired linearity of $c_*$
on $K_p(B,B^+)$.
\end{proof}

\begin{remark}\label{refereetotaro}
Assume $w\in K_p(B,B^+)$ is in the $i$th eigenspace. Then, in the
notation of Proposition \ref{referee} and its proof, we have
\begin{align*}
&\left(\gamma_*[r]c_*\right)(w)=\left([r]\gamma_*
c_*\right)(w)=\left([r]c_*\V_c\gamma_*\right)(w)=\\
&\left(c_*[r]\V_c\gamma_*\right)(w)=
\left(r^{ci}c_*\V_c\gamma_*\right)(w)=\left(r^{ci}\gamma_*\right)(w).
\end{align*}
This shows that the operator $c_*$ sends the $i$th eigenspace in
$K_p(\k[M],\k[M]^+)$ to the $(ci)$th eigenspace. As a result, the
nilpotence conjecture would follow if the graded space
$K_p(\k[M],\k[M]^+)$ were concentrated in finitely many degrees
(with respect to the decomposition as in Proposition
\ref{totaro}(a)). These observations were made by the anonymous
referee and, also, by Burt Totaro in his unpublished notes (1997)
on Conjecture \ref{conj}. Our approach to the problem is through
the finiteness homological result Theorem \ref{findim} below. It
is a weaker statement. But the main result of this section,
Proposition \ref{cyclic}, says that it also suffices for Theorem
\ref{yes} when $\k$ is a number field. This consideration raises
the following interesting question: is the group
$K_p(\k[M],\k[M]^+)$ always concentrated in finitely many degrees
when $M$ is normal?
\end{remark}

\subsection{\emph{Almost isomorphisms}}\label{Almost}
Let ${\bf V}(\Bbb Q)$ denote the category of rational vector
spaces and ${\bf V_f}(\Bbb Q)$ denote the Serre subcategory of
finite dimensional vector spaces. A homomorphism in ${\bf V}(\Bbb
Q )$ will be called an \emph{almost isomorphism} if it defines an
isomorphism in ${\bf V}(\Bbb Q)/{\bf V_f}(\Bbb Q)$, and it will be
called an \emph{almost zero homomorphism} if it defines a zero
homomorphism in ${\bf V}(\Bbb Q)/{\bf V_f}(\Bbb Q)$.

Let $\k$ be a number field (i.~e. a finite extension of $\QQ$),
$A=A_0\oplus A_1\oplus A_2\oplus\cdots$ be a not necessarily
commutative graded $\k$-algebra. In the lemma below and its proof
we use the notation introduced in Section \ref{Linear}.

\begin{lemma}\label{goodwillie}
For any nonnegative integer $p$ and a graded ideal $I\subset A^+$,
such that $\dim_\QQ(A/I)<\infty$, the relative group $K_p(A,I)$ is
a rational vector space and the homomorphism $K_p(A,I)\to
K_p(A,A^+)$ is an almost isomorphism.
\end{lemma}

\begin{proof}
Since $K_p(A_0)$ splits off functorially both from $K_p(A)$ and
$K_p(A/I)$ the long exact sequence, associated to the exact
sequence $0\to I\to A\to A/I\to0$, reads as
$$
\cdots\to K_{p+1}\left(A/I,(A/I)^+\right)\to K_p(A,I)\to
K_p\left(A,A^+\right)\to K_p\left(A/I,(A/I)^+\right)\to\cdots
$$
In particular, by Proposition \ref{totaro}(a) $K_p(A,I)$ is a
$\QQ$-vector space. Furthermore, $(A/I)^+\subset A/I$ is a
nilpotent ideal and by Goodwillie's theorem \cite{Go} we have
$$
K_*(A/I,(A/I)^+)\cong HC_{*-1} (A/I,(A/I)^+),
$$
the groups on the right being finite dimensional $\QQ$-vector
spaces.
\end{proof}

For a ring homomorphism $f:B\subset C$ and an ideal $I\subset B$,
which is mapped bijectively onto an ideal of $C$, the birelative
cyclic homology groups $HC_p(B,C:I)$ are defined as the homology
groups of the shifted (by +1) mapping cone of the homomorphism
between the bicomplexes that produce the relative cyclic homology
groups $HC_p(B,I)$ and $HC_p(C,I)$ (see \cite[\S1]{GW}).

Let $\k$ and $A=A_0\oplus A_1\oplus\cdots$ be as above. Define the
bicomplex $CC_*(A,s)$ by the exact sequences
$$
0\to CC_{*+1}(A_{[s]})\to CC_{*+1}(A)\to CC_*(A,s)\to 0
$$
where $CC_*$ refers to the cyclic bicomplex \cite{Lo}. We will use
the notation $HC_p(A,s)=H_p(CC_*(A,s))$.

\begin{lemma}\label{whitehead}
There exist (natural) homomorphisms
$$
HC_p(A_{[s]},A:A_{[s]}^+)\to HC_p(A,s),\quad p,s\in\ZZ_+,
$$
which are almost isomorphisms when $\dim_\QQ(A/A_{[s]}^+)<\infty$
for all $s\in\ZZ_+$.
\end{lemma}

\begin{proof}
The natural homomorphisms $HC_p(A_{[s]},A:A_{[s]}^+)\to HC_p(A,s)$
exist due to the fact that $HC_*(A,s)$ can be thought of as the
homology groups of the mapping cone of $CC_{*+1}(A_{[s]})\to
CC_{*+1}(A)$. In particular, we have the commutative diagram with
exact rows
$$
\xymatrix{& HC_{p+1}(A_{[s]})\ar[r]&
HC_{p+1}(A)\ar[r]&\\
&HC_{p+1}(A_{[s]},A_{[s]})^+\ar[r]\ar[u]&
HC_{p+1}(A,A_{[s]}^+)\ar[r]\ar[u]&\\
\ar[r] & HC_p(A,s) \ar[r]& HC_p(A_{[s]})\ar[r]&
HC_p(A)\\
\ar[r] & HC_p(A_{[s]},A:A_{[s]}^+) \ar[r]\ar[u]
&HC_p(A_{[s]},A_{[s]}^+)\ar[r]\ar[u]& HC_p
(A,A_{[s]}^+)\ar[u]\\
}
$$
in which the middle vertical homomorphism is an almost isomorphism
because all other vertical homomorphisms are so -- a consequence
of the inequalities
\begin{align*}
\dim_\QQ\left(HC_q(A_{[s]}/A_{[s]}^+)\right),
\dim_\QQ\left(HC_q(A/A_{[s]}^+)\right)<\infty,\quad q\in\ZZ_+.
\end{align*}
\end{proof}

\subsection{\emph{Finite dimensionality
criterion}}\label{fdimensionality}
For the rest of Section \ref{homology} we fix:
\begin{itemize}
\item a number field $\k$,
\item a pyramidal extension of monoids
$M\subset N$, \item a rational cone $D\subset\RR_+M_*$
($=\conv(M_*)$), \item a rational point
$v\in\Phi(N)\setminus\Phi(M)$, \item a nonzero element
$t\in(\RR_+v)\cap\ZZ^n\cong\ZZ_+$.
\end{itemize}

Let $\Lambda=\Lambda(D,t)\subset M_{2\times2}(\k[N])$ denote the
corresponding subring, as in Section \ref{newdescent}.

Fix a grading $\k[N]=\k\oplus B_1\oplus B_2\oplus\cdots$ such that
the elements of $N$ are homogeneous (Section \ref{Monoids}).
Consider the induced gradings
$\Lambda=\Lambda_0\oplus\Lambda_1\oplus\Lambda_2\oplus\cdots$ and
$M_{2\times2}(\k[M])=M_0\oplus M_1\oplus M_2\oplus\cdots$.
Clearly,
$$
\Lambda_0=
\begin{pmatrix}
\k&0\\
\k&\k
\end{pmatrix},\quad
\dim_{\QQ}\Lambda_j<\infty\quad\text{for all}\quad j\ge0,
\quad\text{and}\quad M_0= M_{2\times2}(\k).
$$
Moreover, we have

\begin{equation}\label{assmpt}
\Lambda_s\subset M_s\quad\text{for}\quad s\gg0.
\end{equation}
(See Remark \ref{comments}(a).)

We also fix an arbitrary rational cone $D'$ such that $D\subset
\inte(D')\cup\{O\}$ and $D'\subset\inte(\RR_+M)\cup\{O\}$. We have
the associated monoid $L'=D'\cap\ZZ^n$ and the ring
$$
\Lambda'=
\begin{pmatrix}
\k[L']&\k[L']\cap t^{-1}\k[L']\\
\k[L']+t\k[L']& \k[L']
\end{pmatrix}\subset M_{2\times 2}(\k[N]).
$$
The bigger ring $\Lambda\subset\Lambda'$ also carries the induced
graded structure
$\Lambda'=\Lambda'_0\oplus\Lambda'_1\oplus\Lambda'_2\oplus\cdots$.
Clearly, $\Lambda'_0=\Lambda_0$ and, as in (\ref{assmpt}), we have

\begin{equation}\label{assmpt'}
\Lambda'_s\subset M_s\quad\text{for}\quad s\gg0.
\end{equation}

For arbitrary nonnegative integers $p$ and $s$ we have the
commutative square of group homomorphisms between relative
$K$-groups
\begin{equation}\label{square}
\xymatrix{
K_p(\Lambda_{[s]},(\Lambda_{[s]})^+)\ar[r]^{\theta_s}\ar[d]&
K_p(\Lambda,\Lambda^+)\ar[d]\\
K_p(\Lambda'_{[s]},(\Lambda'_{[s]})^+)\ar[r]_{\theta'_s}&
K_p(\Lambda',(\Lambda')^+), }
\end{equation}
yielding the homomorphisms of rational vector spaces (Lemma
\ref{goodwillie})
$$
\alpha_s:\Coker(\theta_s)\to\Coker(\theta'_s),\quad s\in\ZZ_+.
$$

\begin{proposition}\label{findi}
Conjecture \ref{conj} is true for monoid $\k$-algebras if the
homomorphisms $\alpha_s$ are almost zero for all
$s\gg0$.\footnote{That is, the homomorphisms $\alpha_s$ are almost
zero for arbitrarily fixed $M$, $N$, $D$, $t$ and $p$ as above and
$s\gg0$, depending on these objects.}
\end{proposition}

\begin{proof} Fix a natural number $s$ big enough to satisfy the
following conditions:
\begin{itemize}
\item $\dim_\QQ\I(\alpha_s)<\infty$, \item $\Lambda'_q\subset M_q$
for all $q\ge s$ (see (\ref{assmpt'})).
\end{itemize}

Let a sequence ${\bf c}=(c_1,c_2,\ldots)$ and an element $x\in
K_p(\k[N])$ be such that
$$
(c_1\cdots c_j)_*(x)\in\I\left(K_p(\Lambda)\to K_p(M_{2\times 2
}(\k[N])=K_p(\k[N])\right)
$$
for all $j\gg0$. If we show that
$$
(c_1\cdots c_j)_*(x)\in\I\bigl(K_p(\k[M])\to K_p(\k[N])\bigr)
$$
for all $j\gg0$ then, since the pyramidal extension $M\subset N$
was fixed arbitrarily, Conjecture \ref{conj} follows for monoid
$\k$-algebras by Proposition \ref{descent}.

Replacing $x$ by $x-x(0)$ we can assume $x\in K_p(\k[N],\k[N]^+)$.
Then the elements $(c_1\cdots c_j)_*(x)$, $j\gg0$, are in the
image of the map $K_p(\Lambda,\Lambda^+)\to K_p(\k[N],\k[N]^+)$.

By Proposition \ref{totaro}(a) we have the $\k$-linear
decompositions:
\begin{align*}
&K_p(\k[N],\k[N]^+)=T_1\oplus T_2\oplus\cdots\\
&K_p(\Lambda,\Lambda_+)=U_1\oplus U_2\oplus\cdots\\
&K_p(\Lambda_{[s]}',(\Lambda_{[s]}')^+)=V_{s,1}
\oplus V_{s,2}\oplus\cdots\\
&K_p(\Lambda',(\Lambda')^+)=W_1\oplus W_2\oplus\cdots
\end{align*}
which are functorial with respect to graded $\k$-algebra
homomorphisms.

By the very definition of the endomorphisms $c_*:K_p(\k[N])\to
K_p(\k[N])$ (and the condition $x(0)=0$) we have
$$
c_*(x)\in\I\bigl(K_p(\k[N]_{[c]},(\k[N]_{[c]})^+)\to
K_p(\k[N],\k[N]^+)\bigr),\quad c\in\NN.
$$
Since $\bigoplus_i U_i\to\bigoplus_i T_i$ is a graded
homomorphism, we see from Proposition \ref{totaro}(b) that
\begin{equation}\label{8}
(c_1\cdots c_j)_*(x)\in\I\bigl(0\oplus\cdots\oplus0\oplus
U_{c_1\cdots c_j}\oplus U_{c_1\cdots c_j+1}\oplus\cdots\to
K_p(\k[N],\k[N]^+)\bigr)
\end{equation}
for all indices $j$.

Since the homomorphisms $\bigoplus_i U_i\to\bigoplus W_i$ and
$\bigoplus_i V_{s,i}\to\bigoplus_i W_i$ are graded, by the choice
of $s$ we have $\I(U_i)\subset\I(V_{s,i})$ in $W_i$ for $i\gg0$.
The latter inclusions, together with (\ref{8}), imply
$$
(c_1\cdots
c_j)_*(x)\in\I\bigl(K_p(\Lambda'_{[s]},(\Lambda'_{[s]})^+)\to
K_p(\k[N],\k[N]^+)\bigr),\quad j\gg0.
$$
We are done because $\Lambda'_{[s]}\subset M_{2\times2}(\k[M])$.
\end{proof}

\subsection{\emph{Hochschild homological interpretation}}\label{cyclic}
We define the complexes $C_*(\Lambda,s)$ and $C_*(\Lambda',s)$ by
the exact sequences
\begin{align*}
0\to C_{*+1}(\Lambda_{[s]})\to C_{*+1}(\Lambda)\to C_*(\Lambda,s)\to0,\\
0\to C_{*+1}(\Lambda'_{[s]})\to C_{*+1}(\Lambda')\to
C_*(\Lambda',s)\to0,
\end{align*}
where $C_*(-)$ refers to the Hochschild complex \cite{Lo}.
(Notation as in the previous subsections.) Let
$HH_p(\Lambda,s)=H_p(C_*(\Lambda,s))$ for all $p$ and similarly
for $H_p(C_*(\Lambda',s))$.

\begin{proposition}\label{hochcrit}
In order to prove Conjecture \ref{conj} for monoid $\k$-algebras
it is enough to show that the homomorphisms $HH_p(\Lambda,s)\to
HH_p(\Lambda',s)$ are almost zero for all $p,s\in\ZZ_+$.
\end{proposition}

\noindent\emph{Notice.} In the statement of Proposition
\ref{hochcrit} we could only restrict to $s\gg0$ (in the spirit of
Proposition \ref{findi}), but in Section \ref{finite} the claim
will be proved for all $s\in\ZZ_+$.

\begin{proof}
We have the chain of isomorphisms in the category ${\bf
V}(\QQ)/{\bf V_f}(\QQ)$:
$$
HC_{p-1}(\Lambda,s)\cong
HC_{p-1}(\Lambda_{[s]},\Lambda:\Lambda_{[s]}^+)\cong
K_p(\Lambda_{[s]},\Lambda:\Lambda_{[s]}^+)\cong\Coker(\theta_s)
$$
where the first and the third isomorphisms follow respectively
from Lemmas \ref{whitehead} and \ref{goodwillie} -- it is clear
that the finite dimensionality conditions in those lemmas are
satisfied, and the middle isomorphism follows from Cortinas' proof
of the \emph{KABI-conjecture} \cite[Corollary 0.2]{Cor}. Then, by
Proposition \ref{findi}, Conjecture \ref{conj} follows for monoid
$\k$-algebras if the homomorphisms $HC_p(\Lambda,s)\to
HC_p(\Lambda',s)$ are almost zero for all $p$ and $s$. (Notation
as in Lemma \ref{whitehead}.) Finally, the latter condition
follows from the similar condition on the Hochschild homology
groups by the exact sequences
$$
0\to HC_{p-1}(\Lambda,s)\to HH_p(\Lambda,s)\to
HC_p(\Lambda,s)\to0,
$$
resulting from Conne's exact sequences for the graded
$\k$-algebras $\Lambda$ and $\Lambda_{[s]}$ \cite[Theorem
4.1.13]{Lo}.
\end{proof}

\begin{remark}\label{WhyQ}
The analogous action of $\W(\k)$ on cyclic and Hochschild homology
groups for \emph{$\k$-algebras} has been introduced in \cite{DaW}.
However, we need to resort to number fields in Section
\ref{homology} and, what is more important, in Section
\ref{finite} because of the following: in Goodwillie's and
Cortinas' isomorphisms the cyclic homology groups are those of
$\ZZ$-algebras and, unless $\k$ is a number field, the maps of
homology groups we are interested in seem to be $\k$-linear maps
of \emph{infinite rank}.
\end{remark}

\begin{remark}\label{q}
As it becomes clear in Section \ref{essentially}, we only need to
consider the special case $\k=\QQ$. Moreover, for a number field
$\k$ the Hochcshild and cyclic homology of a $\k$-algebra is the
same over $\QQ$ as over $\k$. But our exposition would remain
literally the same if we had restricted to $\k=\QQ$. So we choose
to work with algebras over $\k$ and vector spaces over $\QQ$.
\end{remark}

\section{Setup for the proof in the arithmetic case}\label{setup}

In this preparatory section we setup the notation and prove
several lemmas to be used in the proof of the Hochschild criterion
in Proposition \ref{hochcrit}.

\subsection{\emph{The data}}\label{DATA}
\underline{For this and the next section we fix the following
objects}:
\begin{itemize}
\item[(i)] a pyramidal extension of monoids $M\subset N$,
identifying $\gp(N)$ with $\ZZ^n$, \item[(ii)] a number field $\k$
and a basis $\Bb(\k)\subset \k$ of $\k$ as a rational vector
space, \item[(iii)] a grading on $\k[N]$, making nontrivial
elements of $N$ homogeneous of positive degree (Section
\ref{Monoids}), \item[(iv)] a rational point
$v\in\Phi(N)\setminus\Phi(M)$ and a nonzero element
$t\in(\RR_+v)\cap\ZZ^n$, \item[(v)] two integers $i\ge1$ and
$s\ge2$ and the notation ${\bf n}=2^{2^{i+1}-1}-2$, \item[(vi)] a
sequence of rational cones $D=D_0,D_1,\ldots,D_{\bf n}=D'$ such
that $D_j\subset\inte(D_{j+1})\cap\{O\}$ for $0\le j\le {\bf n}-1$
and $D'\subset\RR_+M_*$, \item[(vii)] a system of natural numbers
$s=\gamma_{\bf n}<\gamma_{{\bf n}-1}<\cdots<\gamma_1<\gamma_0$
satisfying the condition specified below; see (vii)
(\emph{Continued}) after Lemma \ref{divis}.
\end{itemize}

To each cone $D_j$, $0\le j\le {\bf n}$, and the element $t$ we
associate the subring
$$
\Lambda_j=\Lambda(D_j,t)=\begin{pmatrix}
\k[L_j]&\k[L_j]\cap t^{-1}\k[L_j]\\
\k[L_j]+t\k[L_j]& \k[L_j]
\end{pmatrix}\subset M_{2\times2}(\k[N])
$$
where $L_j=D_j\cap\ZZ^n$ (compare with Section \ref{newdescent}).
Clearly, $\Lambda_0\subset\cdots\subset\Lambda_{{\bf n}}$. All
these rings carry the graded structure induced by that of $\k[N]$.

We will also use the notation $D=D_0$, $D'=D_{{\bf n}}$,
$\Lambda=\Lambda_0$ and $\Lambda'=\Lambda_{{\bf n}}$. Thus the
ring extensions $\Lambda_{j-1}\subset\Lambda_j$, $j\in[1,\bf n]$,
are of the same type as $\Lambda\subset\Lambda'$.

For each $0\le j\le {\bf n}$ the following subset of $\Lambda_j$:
\begin{align*}
&{\bigg\{}b
\begin{pmatrix}
m&0\\
0&0
\end{pmatrix},
b
\begin{pmatrix}
m&0\\
0&m
\end{pmatrix}\ :\ b\in\Bb(\k),\ m\in D_j\cap\ZZ^n{\bigg\}}\ \bigcup\\
&{\bigg\{}b
\begin{pmatrix}
0&0\\
m&0
\end{pmatrix}\ :\ b\in\Bb(\k),\ m\in
\bigl(D_j\cup(t+D_j)\bigr)\cap\ZZ^n{\bigg\}}\ \bigcup\\
&{\bigg\{}b
\begin{pmatrix}
0&m\\
0&0
\end{pmatrix}\ :\ b\in\Bb(\k),\ m\in
D_j\cap(-t+D_j)\cap\ZZ^n{\bigg\}}
\end{align*}
is a basis of $\Lambda_j$ as a rational vector space. We denote
this basis by $\Bb(\Lambda_j)$

\begin{lemma}\label{divis}
For every natural number $\gamma'$ there exists a natural number
$\gamma$ such that any system of homogeneous elements
$\lambda_1,\ldots,\lambda_p\in\Lambda$ of degree $>\gamma$
($p\in\NN$) admits representations of the form
$$
\lambda_1=\lambda_1'\lambda',\ldots,\lambda_p=\lambda_p'\lambda'
$$
where $\lambda_1',\ldots,\lambda_p',\lambda'\in\Lambda'$ are
homogeneous elements of degree $>\gamma'$ and $\lambda'$ is a
central non-zerodivisor. Moreover, if
$\lambda_1,\ldots,\lambda_p\in\Bb(\Lambda)$ then the elements
$\lambda_1',\ldots,\lambda_p'$ and $\lambda'$ can be chosen from
$\Bb(\Lambda')$.
\end{lemma}

\begin{proof}
Fix a natural number $\gamma'$ and a nonzero element $m_0\in
D'\cap\ZZ^n$. We have the easily checked  implication: if
$\deg(m)\gg0$ then the condition $m\in D\cap\ZZ^n$ (or
$m\in(t+D)\cap\ZZ^n$, or $m\in(-t+D)\cap\ZZ^n$) implies
$m_0^{-\gamma'-1}m\in D'\cap\ZZ^n$ (respectively,
$m_0^{-\gamma'-1}m\in(t+D')\cap\ZZ^n$, or
$m_0^{-\gamma'-1}m\in(-t+D')\cap\ZZ^n$).

Therefore, there exists a natural number $\gamma$ such that if
$\deg(\lambda_1),\ldots,\deg(\lambda_p)>\gamma$ then
$$
\lambda_1=\lambda_1'\begin{pmatrix}m_0^{\gamma'+1}&0\\0&m_0^{\gamma'+1}
\end{pmatrix},\ldots,\lambda_p=\lambda_p'
\begin{pmatrix}m_0^{\gamma'+1}&0\\0&m_0^{\gamma'+1}\end{pmatrix}
$$
for some homogeneous elements
$\lambda_1',\ldots,\lambda_p'\in\Lambda'$ of degree $>\gamma'$. In
particular, we have the desired representations in which
$$
\lambda'=\begin{pmatrix}m_0^{\gamma'+1}&0\\0&m_0^{\gamma'+1}
\end{pmatrix}
$$

The second part of the lemma is obvious.
\end{proof}

By applying Lemma \ref{divis} to the successive pairs in the
decreasing sequence $\Lambda_{\bf n}\supset\Lambda_{{\bf
n}-1}\supset\cdots\supset\Lambda_0$ we can choose the natural
numbers $\gamma_j$ so that the following condition holds:
\begin{itemize}
\item[(vii)] (\emph{Continued})
for every index $1\le j\le {\bf n}$ and every system of elements
$\lambda_1,\ldots,\lambda_p\in\Bb(\Lambda_{j-1})$ of degree
$>\gamma_{j-1}$ there exist representations of type
$$
\lambda_1=\lambda_1'\lambda',\ldots, \lambda_p=\lambda_p'\lambda'
$$
where $\lambda_1',\ldots,\lambda_p'$ and $\lambda'$ are elements
of $\Bb(\Lambda_j)$ of degree $>\gamma_j$ and $\lambda'$ is a
central non-zerodivisor.
\end{itemize}

\subsection{\emph{Hochschild complexes}}\label{complexes}
\underline{We will always identify} the $i$-th component
$$
C_i(\Lambda_j,s)\subset C_*(\Lambda_j,s),\quad 0\le j\le {\bf n}
$$
with the rational vector subspace
$$
{\bigg\{}\sum_k
\xi_k\lambda_{k0}\otimes\cdots\otimes\lambda_{ki}{\bigg\}}\subset
\Lambda_j^{\otimes(i+1)}=C_i(\Lambda_j)\subset C_*(\Lambda_j)
$$
(the tensor products being considered over $\QQ$) where
$\xi_k\in\QQ$, $\lambda_{k0},\ldots,\lambda_{ki}\in\mathcal
B(\Lambda_j)$ and
\begin{equation}\label{rose}
\{\deg(\lambda_{k0}),\ldots,\deg(\lambda_{ki})\}\cap\{0,1,
\ldots,s-1\} \not=\emptyset.
\end{equation}

\

\noindent\emph{Caution.} The differential of the complex
$C_*(\Lambda_j,s)$ is \emph{not} the restriction of the
differential of the complex $C_*(\Lambda_j)$.

\

The subset
$$
\Bb(\Lambda_j)^{\otimes(i+1)}:=\{\bar\lambda=\lambda_0\otimes\cdots
\otimes\lambda_i\ :\ \lambda_0,\ldots,\lambda_i\in
\Bb(\Lambda_j)\}\subset C_i(\Lambda_j),\quad 0\le j\le {\bf n},
$$
is a basis of the rational vector space $C_i(\Lambda_j)$, i.~e.
every element $z\in C_i(\Lambda_j)$ ($0\le j\le {\bf n}$) may be
expressed uniquely as a rational linear combination of
$\Bb(\Lambda_j)^{\otimes(i+1)}$; we call this expression the
\emph{canonical expansion} of $z$.

Because of (\ref{rose}), the natural homomorphism
$C_i(\Lambda_j,s)\to C_i(\Lambda_{j+1},s)$ $0\le j\le {\bf n}-1$,
is an embedding. Thus our convention also implies that, as a
rational vector space, $C_i(\Lambda_j,s)$ is identified with its
image in $C_i(\Lambda_{j+1},s)$.

Clearly, the grading on $\Lambda_j$ induces that on
$C_i(\Lambda_j)$ so that
$$
\deg(\lambda_0\otimes\cdots\otimes\lambda_i)=\deg(\lambda_0)+\cdots+
\deg(\lambda_i),\quad \lambda_0,\ldots,\lambda_i\in\Lambda_j.
$$

Recall, the $i$th differential  of the complex $C_*(\Lambda_j,
s)$, $(0\le j\le\bf n)$ is by definition
$\partial_i=\sum_{r=0}^i(-1)^r d_r$ where $d_r:C_i(\Lambda_j,s)\to
C_{i-1}(\Lambda_j,s)$, $0\le r\le i$, are the homomorphisms
between the vertical Cokernels of the diagrams
$$
\xymatrix{
C_i(\Lambda_j)\ar[r]_{\Delta_r}&C_{i-1}(\Lambda_j)\\
C_i((\Lambda_j)_{[s]})\ar[u]\ar[r]_{\Delta_r}&
C_{i-1}((\Lambda_j)_{[s]})\ar[u] }
$$
with
\begin{align*}
\Delta_r(a_0\otimes a_1\otimes\cdots\otimes
a_i)=a_0\otimes\cdots\otimes a_ra_{r+1}\otimes\cdots\otimes
a_i,\quad 0\le r\le i-1,\\
\Delta_i(a_0\otimes a_1\otimes\cdots\otimes a_i)=a_ia_0\otimes
a_1\otimes\cdots\otimes a_{i-1}.
\end{align*}

\subsection{\emph{$\delta$-invariant}}\label{delta}
We adopt the following notation: for a natural number $l$, an
integer $q$, and a system of objects $*_0,\ldots,*_l$, enumerated
this way, we let $\langle q\rangle_{l+1}$ denote the remainder
after dividing $q$ by $l+1$ and let $*_q=*_{{\langle
q\rangle}_{l+1}}$.

\

Assume $0\le j\le {\bf n}$. For an element
$\bar\lambda=\lambda_0\otimes\cdots\otimes\lambda_i\in
\Bb(\Lambda_j)^{\otimes(i+1)}$, admitting an index $0\le i'\le i$
such that $\deg(\lambda_{i'})\le\gamma_j$, we let
$\l_j(\bar\lambda)$ and $\r_j(\bar\lambda)$ denote the integers
determined by the conditions:

\begin{itemize}
\item $0\le\l_j(\bar\lambda)\le\r_j(\bar\lambda)$, \item
$\deg(\lambda_{\l_j(\bar\lambda)}),
\deg(\lambda_{\l_j(\bar\lambda)+1}),\ldots,\deg
(\lambda_{\r_j(\bar\lambda)})>\gamma_j$, \item
$\deg(\lambda_{\l_j(\bar\lambda)-1}),
\deg(\bar\lambda_{\r_j(\bar\lambda)+1})\le \gamma_j$, \item
$\l_j(\bar\lambda)$ (equivalently, $\r_j(\bar\lambda)$) is the
smallest possible number satisfying the preceding conditions.
\end{itemize}
If $\deg(\lambda_{i'})>\gamma_j$ for all indices $0\le i'\le i$
then we put $\l_j(\bar\lambda)=0$ and $\r_j(\bar\lambda)=i$. The
numbers $\l_j(\bar\lambda)$ and $\r_j(\bar\lambda)$ are not
defined if
$\deg(\lambda_0),\deg(\lambda_1),\ldots,\deg(\lambda_i)\le\gamma_j$.

We always have $0\le\l_j(\bar\lambda)\le i$ and
$0\le\r_j(\bar\lambda)\le 2i-1$.\footnote{$\r_j(\bar\lambda)=
2i-1$ if and only if
$\deg(\lambda_0),\ldots,\deg(\lambda_{i-2}),\deg(\lambda_i)>\gamma_j$
and $\deg(\lambda_{i-1})\le\gamma_j$.}

We let $\delta_j(\bar\lambda)=\r_j(\bar\lambda)-\l_j(\bar\lambda)$
if the numbers on the right exist. In this situation
$0\le\delta_j(\bar\lambda)\le i$. In the remaining situation
$\deg(\lambda_0),\deg(\lambda_1),\ldots,\deg(\lambda_i)\le\gamma_j$
we put $\delta_j(\bar\lambda)=-1$.

The definition above can be summarized as follows:
$\delta_j(\bar\lambda)$ is the length of the \emph{first cluster
of $\gamma_j$-high elements} in $\bar\lambda$, using the cyclic
enumeration.

\

The notion of $\delta$-invariant can be extended to all elements
of $C_i(\Lambda_j)$ as follows: for a nonzero element $z\in
C_i(\Lambda_j)$, whose canonical expansion is
$z=\sum_k\xi_k\bar\lambda_k$, $\xi_k\in\QQ$, we put
$\delta_j(z)=\min_k(\delta_j(\bar\lambda_k))$ (this may be $-1$).
We also define $\delta_j(0)=-1$.

Because $s\leq\gamma_j$ for all $0\le j\le {\bf n}$ by
(\ref{rose}) we have the implications
\begin{equation}\label{implication}
z\in C_i(\Lambda_j,s)\quad \Longrightarrow\quad
-1\le\delta_j(z)\le i-1\quad (0\le j\le {\bf n}).
\end{equation}

\subsection{\emph{Format}}\label{Format}
Let $l$ be a natural number. We say that two elements
$$
\bar\lambda=\lambda_0\otimes\cdots\otimes\lambda_l,\
\bar\lambda'=\lambda'_0\otimes\cdots\otimes\lambda'_l\in
\Bb(\Lambda_j)^{\otimes(l+1)},\quad 0\le j\le {\bf n},
$$
with the property
$\delta_j(\bar\lambda),\delta_j(\bar\lambda')\ge0$ are \emph{of
the same format} if the following conditions hold:
\begin{itemize}
\item $\l_j(\bar\lambda)=\l_j(\bar\lambda')$, \item
$\r_j(\bar\lambda)=\r_j(\bar\lambda')$, \item
$\lambda_u=\lambda'_u$ for all indices
$u\in\{0,\ldots,l\}\setminus
\{\langle\l_j(\bar\lambda)\rangle_{l+1},
\ldots,\langle\r_j(\bar\lambda)\rangle_{l+1}\}.$
\end{itemize}
The corresponding equivalence class of $\bar\lambda$ will be
denoted by $\mathfrak F(\bar\lambda)$.

\

For any basis element $\bar\lambda\in C_i(\Lambda_j,s)$ with
$\delta_j(\bar\lambda)\ge0$ and any index $\l_j(\bar\lambda)\le
r\le\r_j(\bar\lambda)-1$, such that $d_r(\bar\lambda)\not=0$, the
summands in the canonical expansion of $d_r(\bar\lambda)$ share
the format. (By convention, we think of $C_i(\Lambda_j,s)$ as
certain rational vector subspace of $C_i(\Lambda_j)$; see Section
\ref{complexes}). In fact, the product of two elements of
$\Bb(\Lambda_j)$ is a rational linear combination of elements of
$\Bb(\Lambda_j)$ of \emph{same} degree. So the pattern of degrees
of the tensor factors in the summands of $d_r(\bar\lambda)$ are
same. Therefore, the numbers $\l_j(d_r(\bar\lambda))$ and
$\r_j(d_r(\bar\lambda))$ can be defined in a natural way, and the
fact that the mentioned summands have same format follows from the
observation that the first cluster of $\gamma_j$-high elements in
these summands, in the sense of the cyclic enumeration, comes from
that in $\bar\lambda$ (compare with Lemma \ref{jumpdelta}(c)
below). So we define the format $\mathfrak F(d_r(\bar\lambda))$ to
be that of any summand in the canonical expansion of
$d_r(\bar\lambda)$.

\

\noindent\emph{Notice.} That the canonical expansion of
$d_r(\bar\lambda)$ may involve more than one summand is a
consequence of the fact that the basis $\Bb(\Lambda_j)$ is not a
multiplicatively closed set for, say, the matrix
$\begin{pmatrix}0&0\\0&1\end{pmatrix}=
\begin{pmatrix}0&0\\1&0\end{pmatrix}
\begin{pmatrix}0&1\\0&0\end{pmatrix}
$ is never in $\Bb(\Lambda_j)$. Also, $\mathcal B(\k)$ is usually
far from being multiplicatively closed.

\

Let $-^\tau$ refer to the cyclic operator
$(\lambda_0\otimes\cdots\otimes\lambda_i)^\tau=\lambda_i\otimes
\lambda_0\otimes\cdots\otimes\lambda_{i-1}$.

The proof of the following lemma is straightforward
\begin{lemma}\label{bigstar} Assume we are given:
\begin{itemize}
\item $j\in\{0,\ldots,{\bf n}\}$, \item
$\bar\lambda,\bar\mu\in\Bb^{\otimes(i+1)}(\Lambda_j)\cap
C_i(\Lambda_j,s)$ such that $\delta_j(\bar\lambda)>0$ and
$\delta_j(\bar\mu)>0$, \item
$u\in\{\l_j(\bar\lambda),\l_j(\bar\lambda)+1,\ldots,
\r_j(\bar\lambda)-1\}$,
\item
$v\in\{\l_j(\bar\mu),\l_j(\bar\mu)+1,\ldots,\r_j(\bar\mu)-1\}$.
\end{itemize}
Assume further $d_u(\bar\lambda),d_v(\bar\mu)\not=0$ and
$\mathfrak F\left(d_u(\bar\lambda)\right)=\mathfrak F
\left(d_v(\bar\mu)\right)$. Then one of the following seven
conditions holds:
\begin{enumerate}
\item[(a)]
\begin{align*}
\l_j(\bar\lambda)=\l_j(\bar\mu),\ \r_j(\bar\lambda)=\r_j(\bar\mu),\\
\r_j(\bar\lambda)\le i,\ \mathfrak F(\bar\lambda)=\mathfrak
F(\bar\mu).
\end{align*}
\item[(b)]
\begin{align*}
\l_j(\bar\lambda)=\l_j(\bar\mu),\ \r_j(\bar\lambda)=\r_j(\bar\mu),\\
\r_j(\bar\lambda)\ge i+1,\ u,v\le i,\\
\mathfrak F(\bar\lambda)=\mathfrak F(\bar\mu).
\end{align*}
\item[(c)]
\begin{align*}
\l_j(\bar\lambda)=\l_j(\bar\mu),\ \r_j(\bar\lambda)=\r_j(\bar\mu),\\
\r_j(\bar\lambda)\ge i+1,\ u,v\ge i+1,\\
\mathfrak F(\bar\lambda)=\mathfrak F(\bar\mu).
\end{align*}
\item[(d)]
\begin{align*}
\l_j(\bar\lambda)=\l_j(\bar\mu)-1,\ \r_j(\bar\lambda)=\r_j(\bar\mu)-1,\\
\r_j(\bar\lambda)\ge i+1,\ u\le i,\ i+1\le v,\\
\mathfrak F(\bar\lambda^\tau)=\mathfrak F(\bar\mu).
\end{align*}
\item[(e)]
\begin{align*}
\l_j(\bar\lambda)=\l_j(\bar\mu)+1,\ \r_j(\bar\lambda)=\r_j(\bar\mu)+1,\\
\r_j(\bar\lambda)\ge i+2,\ u\ge i+1,\ v\le i,\\
\mathfrak F(\bar\mu^\tau)=\mathfrak F(\bar\lambda).
\end{align*}
\item[(f)]
\begin{align*}
\l_j(\bar\lambda)=0,\ \l_j(\bar\mu)=i,\ \r_j(\bar\lambda)+i=
\r_j(\bar\mu),\\
u\le i,\ v=i,\\
\mathfrak F(\bar\mu^\tau)=\mathfrak F(\bar\lambda).
\end{align*}
\item[(g)]
\begin{align*}
\l_j(\bar\lambda)=i,\ \l_j(\bar\mu)=0,\ \r_j(\bar\lambda)=
\r_j(\bar\mu)+i,\\
u= i,\ v\le i,\\
\mathfrak F(\bar\lambda^\tau)=\mathfrak F(\bar\mu).
\end{align*}
\end{enumerate}
\end{lemma}

\

Consider a nonempty set of integers
$$
\mathcal S=\{s_1,\ldots,s_{\#(\mathcal S)}\},\quad
s_1<\cdots<s_{\#(\mathcal S)},
$$
representing distinct residue classes modulo $i+1$. For an element
$$
\bar\lambda=\lambda_0\otimes\lambda_1\otimes\cdots\otimes\lambda_l\in
\Bb(\Lambda_j)^{\otimes(i+1)},\quad 0\le j\le {\bf n}
$$
we put
$$
\bar\lambda|_{\mathcal
S}=\lambda_{s_1}\otimes\cdots\otimes\lambda_{s_{\#(\mathcal S)}}
$$
In the special case $\delta_j(\bar\lambda)\ge0$ and $\mathcal
S=\{\l_j(\bar\lambda),\l_j(\bar\lambda)+1,\cdots,\r_j(\bar\lambda)\}$
we will use the notation
$\bar\lambda|_{\delta_j}=\bar\lambda|_{\mathcal S}$.

More generally, assume $\sum_k\xi_k\bar\lambda_k$ is the canonical
expansion of a nonzero element $z\in C_i(\Lambda_j)$. Then we put
$$
z|_{\mathcal S}=\sum_k\xi_k(\bar\lambda_k|_{\mathcal S})\in
C_{\#(\mathcal S)-1}(\Lambda_j).
$$
If $\delta_j(z)\ge0$ and the elements $\bar\lambda_k$ are all of
the same format then we put
$$
z|_{\delta_j}=\sum_k\xi_k(\bar\lambda_k|_{\delta_j})\in
C_{\delta_j(z)}(\Lambda_j).
$$

\begin{lemma}\label{delta0}
Let $m\in\NN$, $0\le j\le {\bf n}$ and $\mathcal S\subset\ZZ$ be a
subset as above. Assume $\bar\lambda_1,\ldots,\bar\lambda_m\in
\Bb(\Lambda_j)^{\otimes(i+1)}$ are not necessarily distinct
elements such that $\bar\lambda_1|_{\mathcal
T}=\cdots=\bar\lambda_m|_{\mathcal T}$ for $\mathcal
T=\{0,\ldots,i\}\setminus\{\langle s\rangle_{i+1}\ :\ s\in\mathcal
S\}$. Then for arbitrary rational numbers $\xi_1,\ldots,\xi_m$ we
have the equivalence
$$
\xi_1\bar\lambda_1+\cdots+\xi_m\bar\lambda_m=0\ \iff\
\xi_1(\bar\lambda_1|_{\mathcal S
})+\cdots+\xi_m(\bar\lambda_m|_{\mathcal S})=0,
$$
the equality on the right hand side being considered in
$C_{\#(\mathcal S)-1}(\Lambda_j)$.

In particular, if $\delta_j(\bar\lambda_1)\ge0$ and $\mathfrak
F(\bar\lambda_1)=\cdots=\mathfrak F (\bar\lambda_k)$ then for
arbitrary rational numbers $\xi_1,\ldots,\xi_m$ we have the
equivalence
$$
\xi_1\bar\lambda_1+\cdots+\xi_m\bar\lambda_m=0\ \iff\
\xi_1(\bar\lambda_1|_{\delta_j})+\cdots+\xi_m
(\bar\lambda_m|_{\delta_j})=0.
$$
\end{lemma}

\begin{proof}
A permutation of the tensor factors yields an automorphism of the
rational vector space $C_i(\Lambda_j)$. This automorphism respects
the basis $\Bb(\Lambda_j)^{\otimes(i+1)}\subset C_i(\Lambda_j)$.
There is, therefore, no loss of generality in assuming that
$\mathcal S=\{0,1,\ldots,m\}$ for some $0\le m\le i$. We can also
assume that $m<i$ for otherwise there is nothing to prove.

The Hochschild complex $C_*(\Lambda_j)$ is isomorphic to the free
algebra of non-com\-mu\-ta\-ti\-ve polynomials over $\QQ$ having
no nonzero constant terms, whose variables are labeled by elements
of $\Bb(\Lambda_j)$. The set $\Bb(\Lambda_j)^{\otimes(i+1)}$
corresponds to the set of monomials of degree $i+1$ (with respect
to the grading where the variables have degree 1). Therefore, the
lemma translates into the following obvious equivalence
$$
(\xi_1A_1+\cdots+\xi_mA_m)B=0\ \iff\ \xi_1A_1+\cdots+\xi_mA_m=0
$$
where $A_1,\ldots,A_m$ and $B$ are the appropriate nonzero
monomials in the mentioned free algebra.

That the second part of the lemma is a special case of the first
part follows from the equivalence:
\begin{align*}
&\delta_j(\bar\lambda_1)\ge0,\qquad\mathfrak
F(\bar\lambda_1)=\cdots=\mathfrak F
(\bar\lambda_k)\qquad\Longleftrightarrow\\
&\bar\lambda_1|_{\mathcal T}=\cdots=\bar\lambda_m|_{\mathcal
T},\qquad\mathcal T=\{0,\ldots,i\}\setminus\{\langle
s\rangle_{i+1}\ :\
 s\in\mathcal S\}
\end{align*}
where $\mathcal
S=\{\l_j(\bar\lambda_1),\l_j(\bar\lambda_1)+1,\cdots,
\r_j(\bar\lambda_1)\}$.
\end{proof}

\subsection{\emph{$\delta$-machines}}\label{machine}
By an \emph{$i$-sequence} we mean an element
$(s_0,s_1,\ldots,s_i)\in\{+,-\}^{i+1}$. We follow our convention
on cyclic enumeration modulo $i+1$, see Subsection \ref{delta}.
For an $i$-sequence $\sigma=(s_0,\ldots,s_i)$ containing ``$-$'' a
\emph{cluster in $\sigma$} means a subsequence
$(s_q,s_{q+1},\ldots,s_r)$, enumerated this way, where
$\{q,q+1,\ldots,r\}$ is a maximal set with the properties $0\le
q\le r\le 2i-1$, $q\le i$, and $s_q,s_{q+1},\ldots,s_r=+$. Also,
we let the constant sequence $(+,\ldots,+)$ be the only cluster in
itself.

For an $i$-sequence $\sigma=(s_0,\ldots,s_i)$, containing both
``$+$'' and ``$-$'', we define the numbers $\l(\sigma)$ and
$\r(\sigma)$ as follows:

\begin{itemize}
\item $0\le\l(\sigma)\le\r(\sigma)$, \item $s_{\l(\sigma)},
s_{\l(\sigma)+1},\ldots,s_{\r(\sigma)}=+$, \item
$s_{\l(\sigma)-1}, s_{\r(\sigma)+1}=-$, \item $\l(\sigma)$
(equivalently, $\r(\sigma)$) is the smallest possible number
satisfying the preceding conditions.
\end{itemize}
(Compare with the definition of $\l_j$ and $\r_j$ in Section
\ref{delta}.)

If $s_0=\cdots=s_i=+$ then we put $\l(\sigma)=0$ and
$\r(\sigma)=i$. The numbers $\l(\sigma)$ and $\r(\sigma)$ are not
defined if $s_0,s_1,\ldots,s_i=-$.

If $\sigma$ has ``$+$'' then its \emph{initial cluster}, i.~e. the
cluster that starts with $+_{\l(\sigma)}$, will be denoted by
$\sigma_{\ini}$. If $\sigma=(-,\ldots,-)$ then $\sigma_{\ini}$ is
not defined.

Let $\delta(\sigma)=\r(\sigma)-\l(\sigma)$ if the numbers on the
right exist. In the remaining case $\sigma=(-,\ldots,-)$ we let
$\delta(\sigma)=-1$.

\

\noindent\emph{Observation.} For the $5$-sequence
$\sigma=(+,+,-,+,-,+)$ we have $\l(\sigma)=\r(\sigma)=3$,
$\sigma_{\ini}=(+_3)$, and $\delta(\sigma)=0$. Deleting $+_3$ from
$\sigma$ we get the $4$-sequence $\sigma'=(+,+,-,-,+)$ with
$\sigma'_{\ini}=(+_4,+_5,+_6)$. The following lemma, whose proof
is straightforward, tells us that we have more control on the
behavior of the initial clusters when we delete ``+'' from a
cluster of length $\ge2$.

\begin{lemma}\label{jumpdelta}
Assume $\sigma$ is an $i$-sequence. Let $\sigma'$ be an
$(i-1)$-sequence obtained from $\sigma$ by deleting one ``$+$''
with index $\not=0,i+1$ from a cluster of length $\ge2$. Then the
following are true:
\begin{itemize}
\item[(a)] $\sigma'_{\ini}$ is not the image of $\sigma_{\ini}$ in
$\sigma'$ if and only if $(+_i,\ldots,+_r)$, is a cluster in
$\sigma$, different from $\sigma_{\ini}$, and the deleted ``$+$''
is $+_i$,
\item[(b)] If a cluster in $\sigma$ contains $+_0$ or $+_{i+1}$
then the image cluster in $\sigma'$ contains correspondingly $+_0$
or $+_i$.
\item[(c)] If ``$+$'' was deleted from $\sigma_{\ini}$ then
$\sigma'_{\ini}$ is the image of $\sigma_{\ini}$.
\end{itemize}
\end{lemma}
\noindent(Here the ``image'' of a cluster in $\sigma$ means the
cluster in $\sigma'$ whose ``$+$''-s come from the mentioned
cluster.)

A \emph{contraction} of $\sigma=(s_0,\ldots,s_i)$ is an element of
$\{+,-\}^{i+1}$ obtained from $\sigma$ in one of the following
ways:
\begin{itemize}
\item deleting one ``$-$'', \item deleting one ``$+$'' with index
$\not=0,q,i+1$ in a cluster $(+_q,+_{q+1},\ldots,+_r)$ different
from $\sigma_{\ini}$, where $q<i$ and $q<r$, \item deleting one
``$+$'' with index $\not=i+1$ in a cluster
$(+_i,+_{i+1},\ldots,+_r)$ different from $\sigma_{\ini}$, where
$i<r$.
\end{itemize}

A \emph{transformation} of an $i$-sequence $\sigma$ with the
property $0\le\delta(\sigma)<i$ is an $i$-sequence $\sigma'$ which
is obtained from $\sigma$ by first contracting $\sigma$ and then
enlarging the corresponding image of $\sigma_{\ini}$ in the
contraction by inserting one ``$+$'' on the right of the image of
$+_{\r(\sigma)}$.

Transformations are not defined for $(-,\ldots,-)$ and, by
convention, the constant sequence $(+,\ldots,+)$ is the only
transformation of itself. The set of all transformations of an
$i$-sequence $\sigma$ will be denoted by $\Trf(\sigma)$.

Consider the relation on $i$-sequences: $\sigma\preceq\sigma'$ if
either $\sigma'$ has more ``$+$''-s than $\sigma$ or $\sigma'$ and
$\sigma$ have same number of ``$+$''-s and $\sigma'$ is obtained
from $\sigma$ by a (maybe empty) series of transformations.

\begin{lemma}\label{order}
$\preceq$ is a partial order on the set of $i$-sequences.
\end{lemma}

\begin{proof}
Since a transformation of an $i$-sequence has at least as many
``$+$''-s as the original sequence, we only need to exclude the
existence of a family of $i$-sequences $\sigma_0,\ldots,\sigma_k$,
$k\ge1$, such that:
\begin{itemize}
\item $\sigma_0,\ldots,\sigma_k$ have same number of ``$+$''-s,
\item every $\sigma_j$ is different from $(+,\ldots,+)$ and
$(-,\ldots,-)$, \item $\sigma_{j+1}$ is a transformation of
$\sigma_j$ for $j=0,\ldots,k$, in the sense of the cyclic
enumeration modulo $k$.
\end{itemize}
Since there are only finitely many $i$-sequences it suffices to
prove that there is no infinite series of $i$-sequences
$(\sigma_0,\sigma_1,\ldots)$, all different from $(+,\ldots,+)$,
such that $\delta(\sigma_0)\ge0$ and for every $j\in\ZZ_+$ the
$i$-sequence $\sigma_{j+1}$ is a transformation of $\sigma_j$ with
same number ``$+$''-s.

Assume to the contrary that such a series
$(\sigma_0,\sigma_1,\ldots)$ exists. Then for every $j\in\ZZ_+$
the $i$-sequence $\sigma_{j+1}$ is obtained from $\sigma_j$ by
deleting ``$+$'' from a cluster $\not=(\sigma_j)_{\ini}$ with at
least $2$ terms and adding ``+'' to $(\sigma_j)_{\ini}$. If for
some $j\in\ZZ_+$ the cluster which got the new ``$+$'' is
$(\sigma_{j+1})_{\ini}$ then its image in $\sigma_{j+2}$ again
gets a new ``$+$'', etc. This process cannot continue infinitely
for the number of ``$+$''-s cannot grow infinitely. This means
that there is an index $j$ such that the cluster that got the
extra ``$+$'' in $\sigma_{j+1}$ is not $(\sigma_{j+1})_{\ini}$. By
Lemma \ref{jumpdelta}(a) this can only happen if
$(+_i,+_{i+1},\ldots,+_r)$ is a cluster in $\sigma_j$ which is
different from $(\sigma_j)_{\ini}$ and $(\sigma_j)_{\ini}$ gets a
new ``$+$'' in $\sigma_{j+1}$. But then the $i$th component of
$\sigma_{j+1}$ is ``$-$'', and it remains so in every
$\sigma_{j'}$ with $j'\ge j$ because the contractions involved in
our transformations only affect ``$+$''-s. Therefore, for
\emph{every} $j'>j$ the cluster $(\sigma_{j'+1})_{\ini}$ is
obtained from $(\sigma_{j'})_{\ini}$ by adding ``$+$'' -- the same
contradiction.
\end{proof}

An \emph{improvement} of an $i$-sequence $(s_0,\ldots,s_i)$ is an
$i$-sequence $(s'_0,\ldots,s'_i)$ such that $s'_q=+$ whenever
$s_q=+$ for $q=0,\ldots,i$. For two systems of $i$-sequences $\SS$
and $\SS'$ we say that $\SS'$ is an \emph{improvement} of
$\mathfrak S$ if every element of $\mathfrak S'$ is an improvement
of some element of $\mathfrak S$. In particular, if $\mathfrak S'$
is an improvement of a subsystem of $\mathfrak S$, then it is also
an improvement of $\mathfrak S$.

\begin{definition}\label{machines}
A \emph{$\delta$-machine of rank $i$} is an infinite series
$\bar\SS=\{\SS_0,\SS_1,\ldots\}$ where each $\SS_j$ is a set of
$i$-sequences and the following conditions hold:
\begin{itemize}
\item[(i)] $\delta(\sigma)\ge0$ for every element
$\sigma\in\SS_0$, \item[(ii)] $(+,\ldots,+)\in\SS_j$ for all
$j\in\ZZ_+$, \item[(iii)] for every $j\in\ZZ_+$ either
$\SS_j=\{(+,\ldots,+)\}$ or there is a \emph{non-empty} subset
$\mathfrak T_j\subset\SS_j\setminus\{(+,\ldots,+)\}$ such that
$\SS_{j+1}$ is an improvement of the system
$$
\left(\bigcup_{\sigma\in\mathfrak T_j}\Trf(\sigma)\right)\
\bigcup\ \left(\SS_j\setminus\mathfrak T_j\right)
$$
\end{itemize}
\end{definition}

It is obvious that the condition (i) in a $\delta$-machine
$\bar\SS=(\SS_0,\SS_1,\ldots)$ transfers to the elements of
$\SS_j$ for all $j\in\ZZ_+$.

\begin{lemma}\label{termination}
For every natural number $i$ and every $\delta$-machine
$\bar\SS=(\SS_0,\SS_1,\ldots)$ of rank $i$ we have
$\SS_j=\{(+,\ldots,+)\}$ whenever $j\ge 2^{2^{i+1}-1}-1$.
\end{lemma}

\begin{proof}
First observe that nothing changes in Definition
\ref{machines}(iii) if instead of the set
$\left(\bigcup_{\sigma\in\mathfrak T_j}\Trf(\sigma)\right) \
\bigcup\ \left(\SS_j\setminus\mathfrak T_j\right)$ we use
$\left(\bigcup_{\sigma\in\mathfrak
T_j}\Trf(\sigma)\cup\{(+,\ldots,+)\}\right) \ \bigcup\
\left(\SS_j\setminus\mathfrak T_j\right)$. Also, by Lemma
\ref{order}, we can augment $\preceq$ to a linear order on
$\{+,-\}^{i+1}$. Denote the augmented order by $\le$. We arrive at
the situation where if $\SS_j\not=\{(+,\ldots,+)\}$ then every
element of $\sigma\in\mathfrak T_j$ is changed in $\SS_{j+1}$ by a
family of $i$-sequences $\tau$ such that $\sigma\le\tau$ for each
$\tau$. Now we are done by the following lemma.
\end{proof}

\begin{lemma}\label{sublemma}
Let $(X,\le)$ be a finite linearly ordered set and $\bar{\mathfrak
X}=(\mathfrak X_0,\mathfrak X_1,\ldots)$ be a sequence of subsets
of $X$ satisfying the following condition for every $j\in\ZZ_+$:
if $\mathfrak X_j\not=\{\max(X)\}$ then $\mathfrak X_{j+1}$ is
obtained from $\mathfrak X_j$ by changing a nonempty set of
elements $x\in\mathfrak X_j$ to subsets $\mathfrak Z_x\subset X$,
such that $x<z$ for every $z\in\mathfrak Z_x$, and leaving the
other elements of $\mathfrak X_j$ unchanged. Then $\mathfrak
X_j=\{\max(X)\}$ for all $j\ge 2^{\#(X)-1}-1$.
\end{lemma}

\begin{proof}
Let $n=\#(X)$. The case $n=1$ corresponds to the constant sequence
$\mathfrak X_j=\{\max(X)\}$, $j\in\ZZ_+$.

Let $n>1$ and assume we have shown that every sequence as in Lemma
\ref{sublemma}, associated to a linearly ordered set $Y$ with
$m<n$ elements, stabilizes in at most $a_m$ steps, i.~e. all
members with indices $\ge a_m$ are $\{\max(Y)\}$. We can assume
$0=a_1\le a_2\le\cdots\le a_{n-1}$.

If $\min(X)\not\in\mathfrak X_0$ then $\bar{\mathfrak X}$ consists
of subsets of $X\setminus\{\min(X)\}$ and, therefore, stabilizes
in at most $a_{n-1}$ steps.

Let us show that if $\min(X)\in\mathfrak X_k$ then $k\le a_{n-1}$.
It is easily observed that the sequence
$$
\left(\mathfrak X_0\setminus\{\min(X)\},\ldots,\mathfrak
X_k\setminus\{\min(X)\}\right)
$$
satisfies the same condition as the original sequence
$\bar{\mathfrak X}$ for the indices $j=0,\ldots,k-1$. Therefore,
by our assumption, if $k>a_{n-1}$ then $\mathfrak
X_{k-1}\setminus\{\min(X)\}=\{\max(X)\}$. In view of the condition
on $\bar{\mathfrak X}$, the latter equality implies $\mathfrak
X_k\subset X\setminus\{\min(X)\}$ -- a contradiction. Since
$\mathfrak X_j\subset X\setminus\{\min(X)\}$ for all $j>k$ by our
assumption, we get $\mathfrak X_j=\{\max(X)\}$ for all
$j\ge(k+1)+a_{n-1}$. Therefore, $\mathfrak X_j=\{\max(X)\}$ for
all $j>2a_{n-1}+1$. In particular, we can assume $a_n\le
2a_{n-1}+1$. Now the estimate $a_n\le 2^{n-1}-1$ follows from the
equalities $a_1=0=2^{1-1}-1=0$ and $2^j-1=2(2^{j-1}-1)+1,\quad
j\in\NN$.
\end{proof}

\subsection{\emph{From $\lambda$-monomials
to $i$-sequences}}\label{monoseq} To any element $z\in
C_i(\Lambda_j,s)$, $0\le j\le {\bf n}$ we associate a family of
$i$-sequences $\SS(z,j)$ as follows. If
$\bar\lambda=\lambda_0\otimes\cdots\otimes\lambda_i
\in\Bb^{\otimes(i+1)}(\Lambda_j)\cap C_i(\Lambda_j,s)$ then let
$\sigma(\bar\lambda,j)$ be the $i$-sequence obtained from
$\bar\lambda$ by replacing the tensor factors of degree
$>\gamma_j$ by ``$+$'' and the tensor factors of degree
$\le\gamma_j$ by ``$-$''. If $z=\sum_{\mathfrak
K}\xi_k\bar\lambda_k\not=0$ is the canonical $\QQ$-linear
expansion then let
$$
\SS(z,j)=\big\{\sigma(\bar\lambda_k,j)\big\}_{k\in\mathfrak K} \
\cup\ \{(+,\ldots,+)\}.
$$
Finally, we let $\SS(0,j)=\{(+,\ldots,+)\}$. Because of
(\ref{implication}) we have
\begin{equation}\label{implication2}
z\in C_i(\Lambda_j,s),\ \SS(z,j)=\{(+,\ldots,+)\}\quad
\Longrightarrow\quad z=0,\quad0\le j\le {\bf n}.
\end{equation}

\section{Finite dimensional images in Hochschild homology}\label{finite}

As remarked, we follow the notation introduced in Section
\ref{DATA}.

By Proposition \ref{hochcrit} (and the initial remark in Step 1
below), the following theorem implies Conjecture \ref{conj} for
monoid $\k$-algebras.

\begin{theorem}\label{findim}
$\dim_\QQ\I\bigl(HH_i(\Lambda,s)\to
HH_i(\Lambda',s)\bigr)<\infty$.
\end{theorem}

\begin{proof}
\underline{\bf\emph{Step 1.}} We have
$C_*(\Lambda,0)=C_*(\Lambda,1)=0$. Also,
$C_0(\Lambda,s)=\Lambda/\Lambda_{[s]}$ is a finite dimensional
rational space. Therefore, there is no loss of generality in
assuming that $s\ge2$ and $i\ge1$, as mentioned in condition (v)
in Section \ref{setup}.

We let $Z_i(\Lambda_j,s)$ and $B_i(\Lambda_j,s)$, $0\le j\le {\bf
n}$ denote respectively the rational spaces of cycles and
boundaries in $C_i(\Lambda_j,s)$.

First we show that Theorem \ref{findim} follows from the following

\

\noindent{\bf\em Claim A.} For an arbitrary cycle $z\in
Z_i(\Lambda,s)$ for which $\delta_0(z)\ge0$ there exists a
$\delta$-machine $\bar\SS=\{\SS_0,\SS_1,\ldots\}$ of rank $i$
satisfying the condition: for every $1\le j\le {\bf n}$ the image
of $z$ in $Z_i(\Lambda_j,s)$ is homologous to a cycle $z'\in
Z_i(\Lambda_j,s)$ such that $\SS(z',j)=\SS_j$.

\

Assume $z\in Z_i(\Lambda,s)$ is a nonzero cycle whose standard
expansion is $\sum_k\xi_k\bar\lambda_k$. Assume, further,
$\deg(\bar\lambda_k)>(i+1)\gamma_0$ for all $k$. Then for each $k$
there exists an index $p_k\in\{0,\ldots,i\}$ such that
$\deg(\lambda_{k,p_k})>\gamma_0$. (Here $\lambda_{k,p_k}$ is the
$p_k$th tensor factor of $\bar\lambda_k$.) This means that $0\le
\delta_0(\bar\lambda_k)$ for all $k$, or equivalently $0\le
\delta_0(z)$.  Let $\bar\SS$ be a $\delta$-machine satisfying the
condition in Claim A. Then by Lemma \ref{termination} and
implication (\ref{implication2}) there exists $0\le j\le {\bf n}$
such that $z\in Z_i(\Lambda_j,s)$ is a boundary.

Obviously, $0$ and the nonzero cycles $z\in Z_i(\Lambda,s)$, whose
canonical expansion $z=\sum_k\xi_k\bar\lambda_k$, $\xi_k\in\QQ$
satisfy the condition $\deg(\bar\lambda_k)>(i+1)\gamma_0$ for all
indices $k$, constitute a rational subspace $V\subset
Z_i(\Lambda,s)$. What has been said above shows that the image of
$V$ in $Z_i(\Lambda',s)$ is contained in $B_i(\Lambda',s)$. On the
other hand we have $\dim_{\QQ}(Z_i(\Lambda,s)/V)<\infty$. In
particular, there exists a finite dimensional rational subspace
$W\subset Z_i(\Lambda',s)$ such that the image of $Z_i(\Lambda,s)$
in $Z_i(\Lambda',s)$ is contained in $W+B_i(\Lambda',s)$. This
clearly implies Theorem \ref{findim}.

\

\noindent\underline{\bf\emph{Step 2.}} \emph{From now on until the
end of the proof of Theorem \ref{findim} we fix a nonzero element
$z\in Z_i(\Lambda,s)$ such that $\delta_0(z)\ge0$}.

Assume $z=\sum_k\xi_k\bar\lambda_k$, $\xi_k\in\QQ$, is the
canonical expansion, where we use the notation
$\bar\lambda_k=\lambda_{k0}\otimes\cdots\otimes\lambda_{ki}$. Our
goal is to construct a $\delta$-machine as in Claim A.

Let $z_{\min}$ denote the sum of the summands in the canonical
expansion $z=\sum_k\xi_k\bar\lambda_k$ that involve the
\emph{minimal number} of tensor factors of degree $>\gamma_0$. We
write
$$
z_{\min}=z_{\min}'+z_{\min}''
$$
where
\begin{itemize}
\item $z_{\min}'$ is the sum of those summands
$\xi_k\bar\lambda_k$ in the canonical expansion of $z_{\min}$ for
which $1\le\l_0(\bar\lambda_k)\le\r_0(\bar\lambda_k)\le i$,
\item
$z_{\min}''$ is the sum of those summands $\xi_k\bar\lambda_k$ in
the canonical expansion of $z_{\min}$ for which either
$\l_0(\bar\lambda_k)=0$ or $\r_0(\bar\lambda_k)\ge i+1$.
\end{itemize}

The canonical expansions
$$
z_{\min}'=\sum_{k\in
K'}\xi_k\bar\lambda_k\quad \text{and}\quad z_{\min}''=\sum_{k\in
K''}\xi_k\bar\lambda_k
$$
give rise to the following non-empty set of indices:
$$ \mathfrak K=
\begin{cases}
\{k\in K'\ :\ \delta_0(\bar\lambda_k)=\delta_0(z_{\min}')\}\quad
\text{if}\ z_{\min}'\not=0,\\
\{k\in K''\ :\
\delta_0(\bar\lambda_k)=\delta_0(z_{\min}'')\}\quad\text{if}\
z_{\min}'=0.
\end{cases}
$$

\

\noindent\underline{\bf\emph{Step 3.}} By the condition (vii) on
the sequence of natural numbers $\gamma_0,\ldots,\gamma_{\bf n}$
(see Subsection \ref{DATA}), for all indices $k\in\mathfrak K$ we
can fix representations $\lambda_{k,\r_0(\bar\lambda_k)}=
\lambda_{k,\r_0(\bar\lambda_k)}'\lambda$ where
$\lambda,\lambda'_{k,\r_0(\bar\lambda_k)}\in\Bb(\Lambda_1)$,
$\lambda$ is a central non-zerodivisor in $\Lambda_1$, and
$$
\deg(\lambda),\ \deg\lambda'_{k,\r_0(\bar\lambda_k)}>\gamma_1.
$$
For each index $k\in\mathfrak K$ denote by $\hat\lambda_k$ the
element of $C_{i+1}(\Lambda_1)$, obtained from $\bar\lambda_k$ by
substituting $\lambda'_{k,\r_0(\bar\lambda_k)}\otimes\lambda$ for
the tensor factor $\lambda_{k,\r_0(\bar\lambda_k)}$. Clearly, the
condition (\ref{rose}) in Section \ref{complexes} is satisfied for
$\hat\lambda_k$ as well and, hence, $\hat\lambda_k\in
C_{i+1}(\Lambda_1,s)$.

We define the functions $-^{(k)}:\ZZ\to\ZZ$, $k\in\mathfrak K$, by
$$
r^{(k)}=
\begin{cases}
r\ \text{if}\ \r_0(\bar\lambda_k)\le i,\\
r+1\ \text{else}.
\end{cases}
$$
Notice, if $\r_0(\bar\lambda_k)\ge i+1$ for some $k\in\mathfrak K$
then
\begin{equation}\label{1shift}
\langle r^{(k)}\rangle_{i+2}=
\begin{cases}
r+1\ \text{for}\ 0\le r\le i,\\
\langle r\rangle_{i+1}=r-i-1\ \text{for}\ i+1\le r\le 2i-1.
\end{cases}
\end{equation}
For simplicity of notation we let
$\epsilon_k=(-1)^{\langle\r_0(\bar\lambda_k)^{(k)}\rangle_{i+2}+1}$.

Consider the element
$$
\hat z=\sum_{k\in\mathfrak K}\epsilon_k \xi_k\hat\lambda_k\in
C_{i+1}(\Lambda_1,s).
$$
Then Claim A follows from successive application of the following
claim to the ring extensions $\Lambda_j\subset\Lambda_{j+1}$,
$0\le j\le {\bf n}-1$, with use of condition (vii) in Section
\ref{DATA}:

\

\noindent{\bf\em Claim B.} For the element
$z_1=z+\partial_{i+1}(\hat z)\in Z_i(\Lambda_1,s)$ the
$i$-sequence $\SS(z_1,1)$ is obtained from $\SS(z,0)$ by the
process mentioned in Definition \ref{machines}(iii). (As usual,
$\partial_{i+1}$ is the $(i+1)$st differential in the complex
$C_*(\Lambda_1,s)$.)

\

So it is enough to prove Claim B.

\

\noindent\underline{\bf\emph{Step 4.}} We have
$\partial_{i+1}=\sum_{r=0}^{i+1}(-1)^rd_r$ where
$d_r:C_{i+1}(\Lambda_1,s)\to C_i(\Lambda_1,s)$ are certain
homomorphisms (see Subsection \ref{complexes}). By the definition
of the functions $-^{(k)}$ we get
\begin{equation}\label{kthsummand}
\begin{aligned}
&\sum_{k\in\mathfrak K}\xi_k\bar\lambda_k+\sum_{k\in\mathfrak K}
(-1)^{\langle\r_0(\bar\lambda_k)^{(k)}\rangle_{i+2}}
d_{\r_0(\bar\lambda_k)^{(k)}}\bigl(\epsilon_k
\xi_k\hat\lambda_k\bigr)=\\
&\sum_{k\in\mathfrak K}\xi_k\bar\lambda_k-\sum_{k\in\mathfrak K}
d_{\r_0(\bar\lambda_k)^{(k)}}(\xi_k\hat\lambda_k)=0.
\end{aligned}
\end{equation}

For every index $k\in\mathfrak K$ consider the sets
\begin{align*}
&\mathcal S_k=\{\l_0(\bar\lambda_k)^{(k)},
(\l_0(\bar\lambda_k)+1)^{(k)},\ldots,
(\r_0(\bar\lambda_k)-1)^{(k)}\},\\
&\mathcal T_k=\{\l_0(\bar\lambda_k)^{(k)},
(\l_0(\bar\lambda_k)+1)^{(k)},\ldots,(\r_0(\bar\lambda_k)-1)^{(k)},
(\r_0(\bar\lambda_k))^{(k)}\}.
\end{align*}
(We do not exclude the case $\mathcal S_k=\emptyset$, i.~e.
$\delta_0(\bar\lambda_k)=0$ for some $k$.) Assume we have shown
\begin{equation}\label{0sum}
\sum_{k\in\mathfrak K}\sum_{r\in\mathcal S_k} (-1)^{\langle
r\rangle_{i+2}} d_r\bigl(\epsilon_k\xi_k\hat\lambda_k\bigr)=0.
\end{equation}
(Here the summation is considered over those $k$ for which
$\mathcal S_k\not=\emptyset$.) Then (\ref{kthsummand}) and
(\ref{0sum}) imply that
\begin{align*}
z_1=\sum_{k\not\in\mathfrak
K}\xi_k\bar\lambda_k+\sum_{k\in\mathfrak K}\sum_{r\notin\mathcal
T_k} (-1)^{\langle r\rangle_{i+2}}
d_r\bigl(\epsilon_k\xi_k\hat\lambda_k\bigr).
\end{align*}
Since $\gamma_0>\gamma_1$, the latter equality shows that every
element of $\SS(z_1,1)$ is either an improvement of
$\sigma(\bar\lambda_k,0)$ for some $k\notin\mathfrak K$ or an
improvement of a transformation of $\sigma(\bar\lambda_k,0)$ for
some $k\in\mathfrak K$. Clearly, this implies Claim B. So we need
to prove (\ref{0sum}).

\

\noindent\underline{\bf\emph{Step 5.}} Let $d'_r$,
$r\in\{0,\ldots,i\}$ denote the homomorphisms $C_i(\Lambda,s)\to
C_{i-1}(\Lambda,s)$ for which the $i$th differential
$\partial_i:C_i(\Lambda,s)\to C_{i-1}(\Lambda,s)$ is given by
$\partial_i=\sum_{r=0}^i(-1)^rd'_r$.

In this step we prove the following equality:
\begin{equation}\label{cycle0sum}
\sum_{
{\tiny\begin{matrix}k\in\mathfrak K\\
\l_0(\bar\lambda_k)\le r\le \r_0(\bar\lambda_k)-1
\end{matrix}}}(-1)^{\langle r\rangle_{i+1}}\xi_kd'_r(\bar\lambda_k)=0
\end{equation}
\noindent Later we will use it in the proof of (\ref{0sum}).

First consider the case when $z_{\min}'\not=0$ (notation as in
Step 2).

Let $z'$ be a summand in the canonical expansion of some
$d'_r(\xi_k\bar\lambda_k)$, $k\in\mathfrak K$,
$\l_0(\bar\lambda_k)\le r\le \r_0(\bar\lambda_k)-1$. Then by Lemma
\ref{jumpdelta} we have:
\begin{equation}\label{noname}
1\le\l_0(z')\le\r_0(z')\le i-1\quad
\text{and}\quad\delta(\sigma(z',0))=\delta_0(z')=\delta_0(z_{\min}')-1.
\end{equation}
(Notation as in Section \ref{monoseq}.)

For (\ref{cycle0sum}) it is enough to show
\begin{equation}\label{nname}
\sum_{\tiny{\begin{matrix}
k'\in\mathfrak K\\
\l_0(\bar\lambda_{k'})\le
r\le\r_0(\bar\lambda_{k'})-1\\
\sigma(z',0)=\sigma(d'_r(\bar\lambda_{k'}),0)
\end{matrix}}}(-1)^{\langle r\rangle_{i+1}}\xi_{k'}
d'_r(\bar\lambda_{k'})=0
\end{equation}
because the left hand side of (\ref{cycle0sum}) breaks up into
subsums, each of which is of the form the left hand side of
(\ref{nname}) for appropriate $z'$. Since $z$ is a cycle, for
(\ref{nname}) we only need to show the following implication for a
summand $x$ in the canonical expansion of $z$ and any index $0\le
r\le i$:
\begin{equation}\label{nnname}
\begin{aligned}
\sigma(z',0)=\sigma(d'_r(x,0) \quad \Longrightarrow\quad&
x=\xi_{k'}\bar\lambda_{k'}\ \text{for some}\ k'\in\mathfrak K\\
&\text{and}\ r\in\{\langle\l_0(\bar\lambda_{k'})\rangle_{i+1},
\ldots,\langle\r_0(\bar\lambda_{k'})-1\rangle_{i+1}\}.
\end{aligned}
\end{equation}

\

\noindent(Here $\sigma(d'_r(x),0)$ refers to $\sigma(y,0)$ for any
summand $y$ in the canonical expansion of $d'_r(x)$ -- they are
all same.)

\

Let $x$ and $r$ satisfy the condition in (\ref{nnname}). Since
$\bar\lambda_k$ has the smallest possible number of tensor factors
$\lambda_{kr}$ of degree $>\gamma_0$, it is necessary that $x$ is
a summand of $z_{\min}$ and $r$ satisfies the condition
$\deg(\lambda_{kr}),\deg(\lambda_{kr+1})>\gamma_0$. Now if either
$x$ or $r$ does not satisfy the condition in the conclusion of
(\ref{nnname}) then by Lemma \ref{jumpdelta} one of the following
conditions holds:
\begin{itemize}
\item $1\le\l_0(d'_r(x))\le\r_0(d'_r(x))\le i-1$\ \ \ and\ \ \
$\delta(\sigma(d'_r(x),0))=\delta_0(z_{\min}')$,
\item
either\ \  $0=\l_0(d'_r(x))$\ \  or\ \  $\r_0(d'_r(x))\ge i$,
\end{itemize}
-- a contradiction because of (\ref{noname})\footnote{The notation
$\l_0(d'_r(x))$ and $\r_0(d'_r(x))$ we use here is similar to
$\l_j(d_r(\bar\lambda))$ and $\r_j(d_r(\bar\lambda))$ in Section
\ref{Format}, before Lemma \ref{bigstar}.}.

The case $z'_{\min}=0$ follows by essentially the same argument,
with use of the same Lemma \ref{jumpdelta}. One starts with the
observation that for any summand $z''$ in the canonical expansion
of $d'_r(\xi_k\bar\lambda_k)$, $k\in\mathfrak K$,
$\l_0(\bar\lambda_k)\le r\le \r_0(\bar\lambda_k)-1$ the following
conditions hold:
\begin{itemize}
\item
either\ \  $0=\l_0(z'')$\ \  or\ \  $\r_0(z'')\ge i$,
\item
$\delta(\sigma(z'',0))=
\delta_0(z'')=\delta_0(z_{\min}'')-1=\delta_0(z_{\min})-1$.
\end{itemize}

The equality (\ref{cycle0sum}) has been proved.

\

\noindent\underline{\bf\emph{Step 6.}} The set of all pairs
$(k,r)$, $k\in\mathfrak K$, $\l_0(\bar\lambda_k)\le r\le
\r_0(\bar\lambda_k)-1$ with the property
$d'_r(\bar\lambda_k)\not=0$ breaks up into equivalence classes
defined by the relation
$$
(k,r)\sim(l,v)\iff \mathfrak F(d'_r(\bar\lambda_k))=\mathfrak
F(d'_v(\bar\lambda_l)).
$$
Then (\ref{cycle0sum}) and the definition of the format in Section
\ref{Format} imply
$$
\sum_{(k,r)\in\mathfrak C}(-1)^{\langle
r\rangle_{i+1}}\xi_kd'_r(\bar\lambda_k)=0
$$
for every equivalence class $\mathfrak C$.

\

In particular, for every class $\mathfrak C$ by Lemma \ref{delta0}
we have
\begin{equation}\label{0sum''''}
\sum_{(k,r)\in\mathfrak C}(-1)^{\langle
r\rangle_{i+1}}\xi_k\bigl(d'_r(\bar\lambda_k)|_{\delta_0}\bigr)=0.
\end{equation}

We write
\begin{align*}
&\sum_{k\in\mathfrak K}\sum_{r\in\mathcal S_k} (-1)^{\langle
r\rangle_{i+2}} d_r\bigl(\epsilon_k
\xi_k\hat\lambda_k\bigr)=\\
&\sum_{\mathfrak C}\sum_{(k,r)\in\mathfrak C}(-1)^{\langle
r^{(k)}\rangle_{i+2}}\epsilon_k \xi_kd_{r^{(k)}}(\hat\lambda_k)+
\sum_{{\tiny\begin{matrix}k\in\mathfrak K\\\l_0(\bar\lambda_k)\le
r\le \r_0(\bar\lambda_k)-1\\d'_r(\bar\lambda_k)=0
\end{matrix}}}(-1)^{\langle
r^{(k)}\rangle_{i+2}}\epsilon_k \xi_kd_{r^{(k)}}(\hat\lambda_k).
\end{align*}
Therefore, in order to prove (\ref{0sum}) it is enough to prove
the following two things:
\begin{equation}\label{0sum''}
\sum_{(k,r)\in\mathfrak C}(-1)^{\langle
r^{(k)}\rangle_{i+2}}\epsilon_k \xi_kd_{r^{(k)}}(\hat\lambda_k)=0
\end{equation}
for every equivalence class $\mathfrak C$, and
\begin{equation}\label{0}
d'_r(\bar\lambda_k)=0\quad\Longrightarrow\quad
d_{r^{(k)}}(\hat\lambda_k)=0
\end{equation}
for every pair $(k,r)$ such that $k\in\mathfrak K$ and
$\l_0(\bar\lambda_k)\le r\le \r(\bar\lambda_k)-1$.

\

\noindent\underline{\bf\emph{Step 7.}} Fix an equivalence class
$\mathfrak C$ and consider the set
$$
|\mathfrak C|=\{\l_0(d'_r(\bar\lambda_k))^{(k)},\ldots,
\r_0(d'_r(\bar\lambda_k))^{(k)}\}
$$
where $(k,r)\in\mathfrak C$ is any element ($|\mathfrak C|$ is
independent of the choice of $(k,r)$.) It is easily seen that we
are in the situation of Lemma \ref{delta0} with respect to the set
$|\mathfrak C|$ and the elements $d_{r^{(k)}}(\hat\lambda_k)\in
B_i(\Lambda_1,s)$, $(k,r)\in\mathfrak C$. (Notice, here we need
that $\lambda,\lambda'_{k,\r_0(\bar\lambda_k)}\in\Bb(\Lambda_1)$.)
Therefore, by Lemma \ref{delta0} the equality (\ref{0sum''}) is
equivalent to
\begin{equation}\label{0sum'''}
\sum_{(k,r)\in\mathfrak C}(-1)^{\langle r^{(k)}\rangle_{i+2}}
\epsilon_k\xi_k\bigl(d_{r^{(k)}}(\hat\lambda_k)|_{|\mathfrak
C|}\bigr)=0.
\end{equation}
\noindent(Notation as in Lemma \ref{delta0}.)

For every element $(k,r)\in\mathfrak C$ the element
$d_{r^{(k)}}(\hat\lambda_k)|_{|\mathfrak C|}\in
C_{\delta_0(\bar\lambda_k)-1}(\Lambda_1)$ is obtained from
$d'_r(\bar\lambda_k)|_{\delta_0}\in
C_{\delta_0(\bar\lambda_k)-1}(\Lambda)\subset
C_{\delta_0(\bar\lambda_k)-1}(\Lambda_1)$ by substituting
$\lambda'_{k,\r_0(\bar\lambda_k)}$ for the common last tensor
factor $\lambda_{k,\r_0(\bar\lambda_k)}$ in the summands of the
canonical expansion of $d'_r(\bar\lambda_k)|_{\delta_0}$.

The rational vector space
$C_{\delta_0(\bar\lambda_k)-1}(\Lambda_1)$ carries the right
$\Lambda_1$-module structure induced by that on the last (i.~e.
the right-most) tensor factor of
$$
(\Lambda_1)^{\otimes\delta_0(\bar\lambda_k)}=
C_{\delta_0(\bar\lambda_k)-1}(\Lambda_1).
$$
This is a free $\Lambda_1$-module and the element
$\lambda\in\Lambda_1$ is a non-zerodivisor. Therefore, the
equation (\ref{0sum'''}) is equivalent to the equation obtained by
the multiplication of (\ref{0sum'''}) by $\lambda$ in the sense of
the mentioned right $\Lambda_1$-module structure. We will write
(\ref{0sum'''})$\lambda$ for this new equation.

By definition of $\epsilon_k$ in Step 3 (after (\ref{1shift})), if
the parity of the numbers
$$
\Delta(k,r):=\langle
r^{(k)}\rangle_{i+2}+\langle\r_0(\bar\lambda_k)^{(k)}\rangle_{i+2}
-\langle r\rangle_{i+1}
$$
are the same when $(k,r)$ runs through the equivalence class
$\mathfrak C$ then the left hand side of (\ref{0sum'''})$\lambda$
equals the left hand side of (\ref{0sum''''}) and, hence,
(\ref{0sum'''})$\lambda$ follows.

\

\noindent\underline{\bf\emph{Step 8.}} Consider two elements
$(k,r),(k',r')\in\mathfrak C$. Then, according to Lemma
\ref{bigstar}, there are seven possible cases for comparing
$\Delta(k,r)$ and $\Delta(k',r')$. The corresponding computations,
with use of (\ref{1shift}), look as follows:

\

{\tiny\centerline{
\begin{tabular}
{c||c|c}
cases&$(k,r)$&$(k',r')$\\
\hline\hline
(a)&$\begin{matrix}\\r+\r_0(\bar\lambda_k)-r=\r_0(\bar\lambda_k)\\
\
\end{matrix}$&$\begin{matrix}\\r'+\r_0(\bar\lambda_k)-r'=
\r_0(\bar\lambda_k)\\
\ \end{matrix}$\\
\hline (b)&$\begin{matrix}\\(r+1)+
\langle\r_0(\bar\lambda_k)+1\rangle_{i+2}-r=\\\\
\langle\r_0(\bar\lambda_k)\rangle_{i+1}+1\\ \ \end{matrix}$&
$\begin{matrix}\\
(r'+1)+\langle\r_0(\bar\lambda_k)+1\rangle_{i+2}-r'=\\\\
\langle\r_0(\bar\lambda_k)\rangle_{i+1}+1\\ \ \end{matrix}$\\
\hline (c)&$\begin{matrix}\\\langle
r+1\rangle_{i+2}+\langle\r_0(\bar\lambda_k)+1\rangle_{i+2}-\langle
r\rangle_{i+1}=\\\\ \langle\r_0(\bar\lambda_k) \rangle_{i+1}\\ \
\end{matrix}$&$\begin{matrix}\\\langle
r'+1\rangle_{i+2}+\langle\r_0(\bar\lambda_k)+1\rangle_{i+2}-\langle
r'\rangle_{i+1}=\\\\ \langle\r_0(\bar\lambda_k)\rangle_{i+1}\\ \
\end{matrix}$\\
\hline (d)&$\begin{matrix}\\(r+1)+\langle\r_0(\bar\lambda_k)+
1\rangle_{i+2}-r=\\\\
\langle\r_0(\bar\lambda_k)\rangle_{i+1}+1=\r_0(\bar\lambda_k)-i\\
\
\end{matrix}$&$\begin{matrix}\\\langle
r'+1\rangle_{i+2}+\langle\r_0(\bar\lambda_k)+2\rangle_{i+2}-\langle
r'\rangle_{i+1}=\\\\
\langle\r_0(\bar\lambda_k)+1\rangle_{i+1}=\r_0(\bar\lambda_k)-i\\
\
\end{matrix}$\\
\hline (e)&$\begin{matrix}\\\langle
r+1\rangle_{i+2}+\langle\r_0(\bar\lambda_k)^{(k)}
\rangle_{i+2}-\langle r\rangle_{i+1}=\\\\
\r_0(\bar\lambda_k)-i-1\\ \
\end{matrix}$&$\begin{matrix}\\(r'+1)+
\langle(\r_0(\bar\lambda_k)-1)^{(k)}
\rangle_{i+2}-r'=\\\\\r_0(\bar\lambda_k)-i-1\\
\ \end{matrix}$\\
\hline (f)&$r+\r_0(\bar\lambda_k)-r=\r_0(\bar\lambda_k)$
&$\begin{matrix}\\i+1+\langle\r_0(\bar\lambda_k)+i+1\rangle_{i+2}-i=\\\\
\langle\r_0(\bar\lambda_k)+i\rangle_{i+1}+1=\r_0(\bar\lambda_k)\\
\ \end{matrix}$\\
\hline (g)&$\begin{matrix}\\i+1+\langle\r_0(\bar\lambda_{k'})
+i+1\rangle_{i+2}-i=\\\\
\langle\r_0(\bar\lambda_{k'})
+i\rangle_{i+1}+1=\r_0(\bar\lambda_{k'})\\ \
\end{matrix}$&$r'+\r_0(\bar\lambda_{k'})-r'=\r_0(\bar\lambda_{k'})$
\end{tabular}
} }

\

\noindent In particular, not just the residues modulo 2, but the
numbers $\Delta(k,r)$ themselves are the same when $(k,r)$ runs
through the fixed class $\mathfrak C$. This proves (\ref{0sum'''})
and, thus (\ref{0sum''}).

As for the implication (\ref{0}), we can use in the same way Lemma
\ref{delta0} and the right $\Lambda_1$-module structure on
$C_{\delta_0(\bar\lambda_k)-1}(\Lambda_1)$ to derive the
equivalences:
\begin{align*}
&d_{r^{(k)}}(\hat\lambda_k)=0\ \iff d_{r^{(k)}}(\hat\lambda_k)
|_{|\bar\lambda_k|}=0\ \iff\\
&d'_r(\bar\lambda_k)|_{\delta_0}=0\ \iff\ d'_r(\bar\lambda_k)=0,
\end{align*}
where $|\bar\lambda_k|=\{\l_0(d'_r(\bar\lambda_k))^{(k)},\ldots,
\r_0(d'_r(\bar\lambda_k))^{(k)}\}$ and $k\in\mathfrak K$.
\end{proof}

\section{Fields of characteristic 0}\label{essentially}

Here we derive Theorem \ref{yes} from Theorem \ref{findim}. The
latter verifies Conjecture \ref{conj} for number fields. Our
notation is the same as in Theorem \ref{yes}.

\

\noindent\underline{\bf\em Step 1.} Since $K$-groups commute with
filtered colimits and our monoids are filtered unions of affine
positive monoids there is no loss of generality in assuming that
$M$ is an affine positive monoid.

In this step  we prove Conjecture \ref{conj} for a polynomial
coefficient ring $R=\k[\ZZ_+^d]$ where $d\in\NN$ and $\k$ is an
arbitrary number field.

Fix a sequence ${\bf c}=(c_1,c_2,\ldots)$ and an element $x\in
K_p(\k[\ZZ^d\oplus M])$. Considering $x-x(0)$ (with respect to the
natural augmentation $\k[\ZZ^d\oplus M]\to \k$, see Section
\ref{newdescent}) we can assume $x(0)=0$. We already know that
there exists a natural number $j_0$ such that $(c_1\cdots
c_j)^*(x)=0$ for all $j>j_0$, where for a natural number $c$ we
let $c^*$ denote the endomorphism of $K_p(\k[\ZZ^d\oplus M])$
induced by the $\k$-algebra endomorphism $\k[\ZZ^d\oplus M]\to
\k[\ZZ^d\oplus M]$, $l\mapsto l^c$, $l\in\ZZ^d\oplus M$. For each
natural number $c$ we have $c^*=c^\bullet c_*$, where $c_*$ is as
in Theorem \ref{yes} and $c^\bullet:K_p(\k[\ZZ^d\oplus M])\to
K_p(\k[\ZZ^d\oplus M])$ is induced by the $\k$-algebra
endomorphism $\k[\ZZ^d\oplus M]\to \k[\ZZ^d\oplus M]$, $z\mapsto
z^c$ for $z\in\ZZ^d$ and $m\mapsto m$ for $m\in M$. Since
$c^\bullet$ makes $\k[\ZZ^d\oplus M]$ a free module of rank $d^c$
over itself the transfer map for the $K_p$-groups shows
$d^{c_1\cdots c_j}\cdot(c_1\cdots c_j)_*(x)=0$ for $j>j_0$. But
$K_p(\k[\ZZ^d\oplus M])/K_p(\k)$ is a $\QQ$-vector space by
Proposition \ref{totaro}(a) and, hence, $(c_1\cdots c_j)_*(x)=0$
for $j>j_0$, as desired.

\

\noindent\underline{\bf\em Step 2.} Now we prove the nilpotence
conjecture for coefficient rings of the type $R=S^{-1}\k[\ZZ_+^d]$
where $d\in\NN$, $\k$ is a number field, and $S\subset\k[\ZZ_+^d]$
is a multiplicative subset. To this end we need two general facts.

Recall, a ring $A$ is called $K_p$-\emph{homotopy invariant} if
$K_p(A)=K_p(A[\ZZ_+^d])$ for all $d\in\ZZ_+$.

\begin{proposition}[{\cite[Corollary 1.9]{V1}}]\label{vorst}
If $A$ is a $K_p$-homotopy invariant ring then so is its
localization $S^{-1}A$ with respect to any multiplicative subset
$S\subset A$ of non-zerodivisors.
\end{proposition}
\noindent(As Thomason points out in \cite{TT}, the condition that
$S$ consists of non-zerodivisors is superfluous.)

The second result is a variation of the s.~c. \emph{Swan-Weibel
homotopy trick} \cite{A}. We will apply it to certain submonoids
$H\subset\QQ_+$.

\begin{proposition}\label{swan}
Let $H$ be a (not necessarily finitely generated) monoid without
nontrivial units and $A=\bigoplus_H A_h$ be an $H$-graded ring.
Then for any functor $F$ from rings to abelian groups the
following implication holds
$$
F(A)=F(A[H])\ \Longrightarrow\ F(A_1)=F(A).
$$
\end{proposition}
\noindent (Here $1$ denotes the neutral element of $H$, that is
the unit of $A$.)

\begin{proof}
Consider the ring homomorphism $\iota:A\to A[H]$,
$\sum_ka_{h_k}\mapsto h_ka_{h_k}$. Its composite with $A[H]\to A$,
$h\mapsto 1\in A$, is the identity map on $A$. Therefore,
$F(\iota)$ is a group monomorphism. But $F(\pi)$ is an
isomorphism, where $\pi:A[H]\to A$ is the augmentation induced by
$H\setminus\{1\}\to 0\in A$. We see that $F(\pi\iota):F(A)\to
F(A_1)$ is a monomorphism. But it is also an epimorphism because
$A\to A_1$ is a split ring epimorphism.
\end{proof}

For a monoid $H$ and a sequence ${\bf c}=(c_1,c_2,\ldots)$ we let
$H^{\bf c}$ denote the filtered union
$$
\bigcup_{j=1}^\infty\left\{\frac h{c_1\cdots c_j}\ :\ h\in
H\right\}\subset\QQ\otimes H.
$$
Then Conjecture \ref{conj} holds for a coefficient ring $R'$ and a
monoid $M'$ without nontrivial units if and only if
$K_p(R')=K_p(R'[(M')^{\bf c}])$ for every sequence ${\bf
c}=(c_1,c_2,\ldots)$.

Fix a sequence ${\bf c}=(c_1,c_2,\ldots)$. By Step 1, the ring
$\k[\ZZ_+^d][M^{\bf c}]$ is $K_p$-homotopy invariant. Then by
Proposition \ref{vorst}, $R[M^{\bf c}]$ is $K_p$-homotopy
invariant as well. By the filtered limit argument we get
\begin{equation}\label{zbfc}
K_p(R[M^{\bf c}])=K_p(R[M^{\bf c}][\ZZ_+^{\bf c}]).
\end{equation}
The ring $R[M]$ admits a grading $R\oplus R_1\oplus
R_2\oplus\cdots$ where the elements of $M$ are homogeneous
(Section \ref{Monoids}). It follows that $R[M^{\bf c}]$ admits a
$\ZZ_+^{\bf c}$-grading whose zero component is $R$. Therefore,
Proposition \ref{swan} applies to (\ref{zbfc}) and we get
$K_p(R)=K_p(R[M^{\bf c}])$.

\

\noindent\underline{\bf\em Step 3.} By Step 2 Conjecture
\ref{conj} is valid for pure transcendental extensions of $\QQ$.
Therefore, the direct limit argument and the following induction
lemma complete the proof of Theorem \ref{yes}.

\begin{lemma}\label{finitext}
The validity of the nilpotence conjecture for $K_p(\k_1[M])$
transfers to $K_p(\k_2[M])$ for any finite extension
$\k_1\subset\k_2$ of fields of characteristic 0.
\end{lemma}

This lemma will be proved in Step 4 below. We will need the fact
that rational $K$-theory of rings satisfies Galois descent -- a
special case of Thomason's \'etale descent for localized versions
of $K$-theory, first proved in \cite{T} in the smooth case and
then extended to the singular case as an application of the
local-to-global technique of \cite{TT}. The argument below follows
closely \cite[Lemma 2.13]{T} and \cite[Proposition 11.10]{TT}. It
is, of course, important that Thomason's higher $K$-groups agree
with those of Quillen in the affine (or, more generally,
quasiprojective) case \cite[Theorem 7.6]{TT}.

\begin{lemma}\label{galois}
Let $A\subset B$ be a finite Galois extension of noetherian rings
with Galois group $G$. Then
$K_p(A)\otimes\QQ=H^0(G,K_p(B)\otimes\QQ)$.
\end{lemma}

\begin{proof}
The push-out diagram of rings
\begin{equation}\label{1}
\xymatrix{B\ar[r]^{\iota_1}&B\otimes_AB\\
A\ar[r]_{\iota}\ar[u]^\iota&B\ar[u]_{\iota_2}}
\end{equation}
defines a pull-back diagram of the corresponding affine schemes in
the category of all schemes. The latter diagram satisfies
\emph{all} the conditions of Proposition 3.18 in \cite{TT} (the
\emph{base change formula} for singular schemes). Therefore, by
the mentioned proposition we have the equality of endomorphisms
\begin{equation}\label{2}
\iota_\star\iota^\star=(\iota_2)^\star(\iota_1)_\star :K_p(B)\to
K_p(B)
\end{equation}
where $-_\star$ refers to the functorial homomorphism and
$-^\star$ refers to the corresponding transfer map (contrary to
the scheme-theoretical notation in \cite{Q1}\cite{T}\cite{TT}).
Galois theory identifies the diagram (\ref{1}) with the diagram
\begin{equation}\label{3}
\xymatrix{B\ar[r]^\Delta&B^n\\
A\ar[r]_{\iota}\ar[u]^{\iota}&B\ar[u]_{\Delta_G}}
\end{equation}
where $n=\#(G)$, $\Delta$ is the diagonal embedding, and
$\bigl(\Delta_G(b)\bigr)_g=g(b)$. By the elementary properties of
$K$-groups \cite[\S2]{Q1} and the equalities
$(g^{-1})_\star=g^\star$, $g\in G$ the equality (\ref{2}) and the
diagram (\ref{3}) imply
$\iota_\star\iota^\star=\sum_G(g^{-1})_\star=\sum_Gg_\star$. The
other composite $\iota^\star\iota_\star$ is the multiplication by
$[B]\in K_0(A)$ on $K_p(A)$. Since $\rank_{A_\nu}B_\nu=n$ for all
$\nu\in\Spec(A)$ and
\begin{align*}
&K_0(A)_{\text{red}}\cong\{f:\Spec(A)\to\ZZ,\ f\
\text{continuous}\}=
\underbrace{\ZZ\times\ZZ\times
\cdots\times\ZZ}_{\tiny{\#(\text{components
of}\ \Spec(A))}},\\
&[P]\mapsto(\rank_AP:\Spec(A)\to\ZZ),
\end{align*}
(\cite[Ch.9, Proposition 4.4]{Ba}), the class $[B]\otimes\QQ$ is a
unit in $K_0(A)\otimes\QQ$. Therefore, the homomorphisms
$\iota_\star\otimes\QQ$ and the corresponding restriction of
$\left([B]\otimes\QQ\right)^{-1}(\iota^\star\otimes\QQ)$ establish
the desired isomorphism.
\end{proof}

\

\noindent\underline{\bf\em Step 4.} Here we prove Lemma
\ref{finitext}. For clarity we let $c\mapsto c_*/\k$ denote the
induced multiplicative action of $\NN$ on $K_p(\k_1[M])/K_p(\k)$,
$\k$ an arbitrary field.

An embedding of $M$ into a free monoid $\ZZ_+^r$ (see Subsection
\ref{Monoids}) gives rise to compatible graded structures on
$\k_1[M]$ and $\k_2[M]$ such that the monoid elements are
homogeneous. By Proposition \ref{totaro}(a)
$K_p(\k_2[M])/K_p(\k_2)$ carries a $\k_2$-linear structure,
$K_p(\k_1[M])/K_p(\k_1)$ carries a $\k_1$-linear structure, and
the group homomorphism $K_p(\k_1[M])/K_p(\k_1)\to
K_p(\k_2[M])/K_p(\k_2)$ is $\k_1$-linear.

We can also assume that $\k_1\subset\k_2$ is a Galois extension.
In fact, if $\k_1\subset\k_3$ is a finite Galois extension such
that the conjecture is true for monoid $\k_3$-algebras and
$\k_2\subset\k_3$ then the commutative squares
\begin{equation}\label{normal}
\xymatrix{K_p(\k_3[M])/K_p(\k_3)\ar[r]^{c_*/\k_3}&
K_p(\k_3[M])/K_p(\k_3)\\
K_p(\k_2[M])/K_p(\k_2)\ar[r]_{c_*/\k_2}\ar[u]&
K_p(\k_2[M])/K_p(\k_2)\ar[u]
},\qquad c\in\NN
\end{equation}
imply the validity of the nilpotence conjecture for $\k_2$ too
because the vertical homomorphisms in (\ref{normal}) are
monomorphisms. That the mentioned homomorphisms are monomorphisms
is shown as follows: $\k_3[M]$ is a free module over $\k_2[M]$ of
rank $[\k_3:\k_2]$ and, hence, the composite of the functorial map
$K_p(\k_2[M])\to K_p(\k_3[M])$ with the corresponding transfer map
$K_p(\k_3[M])\to K_p(\k_2[M])$ is the multiplication by
$[\k_3:\k_2]$ -- an automorphism on the subgroup
$K_p(\k_2[M])/K_p(\k_2)\subset K_p(\k_2[M])$ which is a rational
vector space.

Let $G$ be the Galois group of the extension $\k_1\subset\k_2$.
Then $\k_1[M]\subset\k_2[M]$ is a Galois extension of rings with
the same Galois group. The action of $G$ on $\W(\k_2)$ shows that
the $\k_2$-vector space $K_p(\k_2[M])/K_p(\k_2)$ is a Galois
$\k_1$-module. Then Lemma \ref{galois} implies
\begin{equation}\label{galdesc}
K_p(\k_2[M])/K_p(\k_2)=K_p(\k_1[M])/K_p(\k_1)\otimes\k_2.
\end{equation}
Since the nilpotence conjecture is valid for $\k_1$, every element
of $K_p(\k_1[M])/K_p(\k_1)$ is annihilated by high iterations of
the multiplicative action of $\NN$. It follows from
(\ref{galdesc}) and the commutative squares of type (\ref{normal})
for the extension $\k_1[M]\subset\k_2[M]$ that the same is true
for the elements of certain generating set of the $\k_2$-vector
space $K_p(\k_2[M])/K_p(\k_2)$. But $c_*/\k_2$ is $\k_2$-linear by
Lemma \ref{referee}. \qed

\end{document}